\documentclass[preprint,sort&compress,12pt]{elsarticle}
\usepackage{amsmath}
\usepackage{mathrsfs}
\usepackage{stmaryrd}
\usepackage{mathrsfs}
\usepackage{bm}
\usepackage{amsmath}
\usepackage{amsfonts}
\usepackage{amssymb}
\usepackage{amsthm}

\journal{\quad }
\newcommand{\mf}{\mathbf}
\newcommand{\mm}{\mathrm}

\begin{document}
\begin{frontmatter}
\title{Nonlinear Instability for Nonhomogeneous \\  Incompressible Viscous Fluids \tnoteref{S}}
 \tnotetext[S]{The research of Fei Jiang
was supported by the NSFC (Grant No. 11101044), the research of Song
Jiang by NSFC (Grant No. 40890154) and the National Basic Research
Program under the Grant 2011CB309705, and the research of Guoxi Ni by NSFC (Grant No. 91130020).}
\author{Fei Jiang\corref{cor1}}
\ead{jiangfei0591@163.com}
 \cortext[cor1]{Corresponding
author: Tel +86-15001201710.}
\author{Song Jiang}
\author{Guoxi Ni}
\address{Institute of Applied Physics and Computational Mathematics, P.O. Box 8009, Beijing 100088, China.}

\begin{abstract}
We investigate the nonlinear instability of a smooth
steady density profile solution of the three-dimensional nonhomogeneous incompressible
Navier-Stokes equations in the presence of a uniform gravitational field,
including a Rayleigh-Taylor steady-state solution
with heavier density with increasing height (referred to the Rayleigh-Taylor instability).
We first analyze the equations obtained from linearization around the
steady density profile solution. Then we construct solutions of the
linearized problem that grow in time in the Sobolev space $H^k$,
thus leading to a global instability result for the linearized
problem. With the help of the constructed unstable solutions and an existence theorem of
classical solutions to the original nonlinear equations, we can then
demonstrate the instability of the nonlinear problem in some sense. Our analysis shows that
the third component of the velocity already induces the instability, this is
different from the previous known results.
\end{abstract}

\begin{keyword}
Nonhomogeneous Navier-Stokes equations, steady density profile, Rayleigh-Taylor instability,
incompressible viscous flows. \MSC[2000] 35Q35\sep  76D03.

\end{keyword}
\end{frontmatter}


\newtheorem{thm}{Theorem}[section]
\newtheorem{lem}{Lemma}[section]
\newtheorem{pro}{Proposition}[section]
\newtheorem{cor}{Corollary}[section]
\newproof{pf}{Proof}
\newdefinition{rem}{Remark}[section]
\newtheorem{definition}{Definition}[section]

\section{Introduction}
\label{Intro} \numberwithin{equation}{section}

This paper is concerned with the nonlinear instability of a smooth
steady density profile solution of the three-dimensional nonhomogeneous incompressible
Navier-Stokes equations in the presence of a uniform gravitational field,
including a Rayleigh-Taylor steady-state solution
with heavier density with increasing height (referred to the Rayleigh-Taylor instability).

 The motion of a nonhomogeneous incompressible viscous fluid in the presence of
a uniform gravitational field in $\mathbb{R}^3$ is governed by the
Navier--Stokes equations:
\begin{equation}\label{0101}\left\{\begin{array}{l}
 \rho_t+{\bf v}\cdot\nabla \rho=0,\\[1mm]
\rho\mathbf{v}_t+\rho {\bf v}\cdot\nabla {\bf v}+\nabla p=\mu\Delta
{\bf v}-\rho{g}e_3,\\[1mm]
\mathrm{div}\mathbf{v}=0,\end{array}\right.\end{equation}
 where the unknowns
$\rho$, $\mathbf{v}$ and $p$ denote the density, the velocity, and
the pressure of the fluid, respectively.
 In the system (\ref{0101}) we have written $\mu>0$ for the
 coefficient of shear viscosity,  $g>0$ for the gravitational
constant, $e_3=(0,0,1)$ for the vertical unit vector, and $-ge_3$
for the gravitational force.

In this paper we consider the problem of the Rayleigh-Taylor (RT) instability, so we assume that
 a smooth steady density profile $\bar{\rho}:=\bar{\rho}({x}_3)\in C^\infty(\mathbb{R})$
exists which satisfies
\begin{eqnarray}\label{0102}
&&\bar{\rho}'\in C_0^{\infty}(\mathbb{R}),\quad \inf_{x_3\in \mathbb{R}}\bar{\rho}>0,\\[0em]
&& \label{0103}\bar{\rho}'(x_3^0)>0\;\;\mbox{ for some point }x_3^0\in
\mathbb{R},
\end{eqnarray}
where $'=d/dx_3$, see Remark \ref{rem0101} on the construction of such $\bar{\rho}$.
 Clearly, such $\bar{\rho}$
with $\mathbf{v}(t,\mathbf{x})\equiv\mathbf{0}$ defines a steady
state to (\ref{0101}), provided
\begin{equation}\label{0104}
\nabla \bar{p}=-\bar{\rho}g e_3,\;\mbox{ i.e., }\;\frac{d\bar{p}}{dx_3}=-\bar{\rho}g .
\end{equation}
\begin{rem}
We point out that by virtue of the condition (\ref{0103}), there is at least
a region in which the steady density solution has larger density
with increasing $x_3$ (height), thus this will lead to the classical
Rayleigh-Taylor instability as shown in Theorem \ref{thm:0101} below.
\end{rem}

 Let the perturbation be
$$ \varrho=\rho -\bar{\rho},\quad \mathbf{u}=\mathbf{v}+\mathbf{0}, \quad q=p-\bar{p},$$
then, $(\varrho ,\mathbf{u},q)$ satisfies the perturbed equations
\begin{equation}\label{0105}\left\{\begin{array}{l}
\varrho_t+{\bf u}\cdot\nabla (\varrho+\bar{\rho})=0, \\[1mm]
(\varrho+\bar{\rho}){\bf u}_t+(\varrho+\bar{\rho}){\bf u}\cdot\nabla
{\bf u}+\nabla q+g \varrho e_3=\mu\Delta{\bf u}, \\[1mm]
\mathrm{div}\mathbf{u}=0.\end{array}\right.  \end{equation}

To complete the statement of the perturbed problem, we specify the initial and boundary
conditions:
\begin{equation}\label{0106}
(\varrho,\  {\bf u} )|_{t=0}=(\varrho_0,\ {\bf u}_0)\quad\mbox{in } \mathbb{R}^3
\end{equation}
and
\begin{equation}\label{0107}
\lim_{|\mf{x}|\rightarrow +\infty}{\bf u}(t,\mathbf{x})={\bf 0}\quad
\mbox{ for any }t>0.
\end{equation}
Moreover, the initial data should satisfy
$\mathrm{div}\mathbf{u}_0=0$.

If we linearize the equations (\ref{0105}) around the steady state
$(\bar{\rho},\mathbf{0})$, then the resulting linearized equations read as
\begin{equation}\label{0108}
\left\{\begin{array}{ll}
 \varrho_t+\bar{\rho}'{u}_3=0, \\[1mm]
  \bar{\rho}\mathbf{u}_t +\nabla q+g\varrho e_3=\mu \Delta \mathbf{u},\\[1mm]
 \mathrm{div}\mathbf{u}=0.
\end{array}\right.\end{equation}

 It has been known for over a century that the steady states $(\bar{\rho},\mathbf{0})$
 to the linearized RT problem (\ref{0106})--(\ref{0108}) with $\mu\geq 0$ is unstable \cite{CSHHS,RLAP},
 i.e., there exists a unstable  solution to (\ref{0106})--(\ref{0108}). Such
 instability to (\ref{0106})--(\ref{0108}) is often called linear RT instability. However, there have
 been only few results on the mathematically rigorous justification of the RT instability for
 (\ref{0105})--(\ref{0107}). In 2003, Hwang and Guo \cite{HHJGY} proved the nonlinear RT
 instability to (\ref{0105})--(\ref{0106}) with boundary condition
 $\mathbf{u}\cdot \mathbf{n}|_{\partial\Omega}=0$ for the two-dimensional inviscid case (i.e. $\mu=0$) where
$\Omega=\{(x_1,x_2)\in \mathbb{R}^2~|~-l<x_2<m\}$ and $\mathbf{n}$ denotes the
 outer normal vector to $\partial\Omega$.
 To our best knowledge, however, it is still open mathematically whether there exists a unstable
solution to the nonlinear RT problem (\ref{0105})--(\ref{0107}) of viscous fluids with
variable density.
The aim of this article is to show rigorously the instability for the nonlinear RT
 problem (\ref{0105})--(\ref{0107}) in some proper sense. The main result read as follows.
\begin{thm}\label{thm:0101}
 Let the steady density profile $\bar{\rho}$ satisfy (\ref{0102})--(\ref{0103}). Then, the
steady state $(\bar{\rho},\mathbf{0})$ of (\ref{0105})--(\ref{0107}) is unstable under the Lipschitz
structure, that is, for any $s\geq 2$, $\delta>0$, $K>0$, and $F$ satisfying
\begin{equation}\label{0109} F(y)\leq Ky\quad\mbox{ for any }y\in [0,\infty),\end{equation}
 there exist a constant $i_0:=i_0(s)>0$ and smooth initial data
\begin{equation*}\label{0110}(\varrho_0,\mathbf{u}_0)\in
(H^\infty(\mathbb{R}^3))^4\mbox{ with }
\|(\varrho_0,\mathbf{u}_0)\|_{H^s(\mathbb{R}^3)}<
\delta,\end{equation*} but the unique classical solution
$(\varrho,\mathbf{u})$ of (\ref{0105})--(\ref{0107}), emanating from
the initial data $(\varrho_0,\mathbf{u}_0)$,  satisfies
\begin{equation}\label{0111}\|{u}_3(t_K)\|_{L^2(\mathbb{R}^3)}>
F(\|(\varrho_0,\mathbf{u}_0)\|_{H^s(\mathbb{R}^3)})\mbox{ for some
}t_K\in
\left(0,\frac{2}{\Lambda}{\mathrm{ln}\frac{2K}{i_0}}\right]\subset
(0,T^{\max}),\end{equation} where the constant $\Lambda$ is given by
(\ref{0267}), $H^\infty(\mathbb{R}^3)=\cap_{k=1}^\infty H^k(\mathbb{R}^3)$,
and $T^{\max}$ denotes the maximal time of existence of the solution $(\varrho,\mathbf{u})$.
\end{thm}
\begin{rem} 
It should be noted that we can not obtain the same instability
result for the domain
$\Omega^{m}_{l}:=\{~\mf{x}\in\mathbb{R}^3~|-l<x_3<m\}$ in place of
$\mathbb{R}^3$, due to lack of an existence result of the classical
solution to the RT problem (\ref{0105}) in the domain
$\Omega^{m}_{l}$. We also mention that Theorem \ref{thm:0101} still
holds if we define
$\|(\varrho_0,\mathbf{u}_0)\|_{H^s(\mathbb{R}^3)}^2:=\|\varrho_0\|_{H^{s-1}(\mathbb{R}^3)}^2
+\|\mathbf{u}_0\|_{H^s(\mathbb{R}^3)}^2$.
\end{rem}
\begin{rem} Our result shows that the problem (\ref{0105})--(\ref{0107}) does not
possess the following stability structure:
\begin{equation}\label{0114}\exists \mbox{ a constant }
C>0,\mbox{ such that }\sup_{t\in (0,
T]}\|{u}_3(t)\|_{L^2(\mathbb{R}^3)}\leq
C\|(\varrho_0,\mf{u}_0)\|_{H^{s}(\mathbb{R}^3)}\;\;\mbox{ for any }T>0,
\end{equation}
which should be quite general and reasonable for a global stability theory.
Notice that $s\geq 2$ in Theorem \ref{thm:0101} is arbitrary. Thus, even if the initial data of the
(\ref{0105})--(\ref{0107}) are smooth and small, the failure of the stability structure
(\ref{0114}) means that it is not possible to use (\ref{0114}) to control
the norm of $\|\mathbf{u}(t)\|_{L^2(\mathbb{R}^3)}$ for long time.
\end{rem}
\begin{rem}\label{rem0101} Here we give an example of a steady density profile
solution satisfying the conditions of Theorem \ref{thm:0101}. Assume
\begin{equation*}\label{}\bar{\rho}^a=\left\{\begin{array}{ll}
                        \tilde{\rho}^{h} & \mbox{ for }x_3\geq 1, \\[1mm]
                        ({\tilde{\rho}^{h}+\tilde{\rho}^{l}})/2 &\mbox{ for } x_3\in (-1,1), \\[1mm]
                         \tilde{\rho}^{l} & \mbox{ for }x_3\leq  -1,
                                        \end{array}\right.
\end{equation*} and $0<\tilde{\rho}^l<\tilde{\rho}^h<+\infty$, then $\bar{\rho}:=S_\varepsilon(\bar{\rho}^a)$
and $\bar{p}:=g\int_{x_3}^0\bar{\rho}(s)\mathrm{d}s$ satisfy
(\ref{0102})--(\ref{0104}), where $S_\varepsilon$ is a standard mollifier operator.
\end{rem}
\begin{rem} The constant $\Lambda$ (in  (\ref{0267})) is often called  maximal linear growth rate.
By virtue of (\ref{0267}) and (\ref{n0229n}), $\Lambda<\infty$, and  $\Lambda\to 0$
if $g\|\bar{\rho}'/\bar{\rho}\|_{L^\infty(\mathbb{R})}\to 0$ or $\mu\to\infty$. In contrast,
$\Lambda=\infty$ in the corresponding inviscid case. This clearly shows that
the viscosity plays an stabilizing role in the linear RT instability.
\end{rem}

 The proof of Theorem \ref{thm:0101}, inspired by \cite{GYTI1,GYTI2}, is divided into four steps:
 (i) First we notice that the coefficients in the linearized equations (\ref{0108})
 depend only on the vertical variable $x_3\in \mathbb{R}^3$, this allows us to seek ``normal mode"
solutions  by taking the horizontal Fourier transform of the
equations and assuming the solutions grow exponentially in time by
the factor $e^{\lambda(|\xi|)t}$, where $\xi\in \mathbb{R}^2$ is the
horizontal spatial frequency and $\lambda(|\xi|)>0$. This reduces
the equations to a system of ordinary differential equations (ODEs) defined on $\mathbb{R}$ with
$\lambda(|\xi|)>0$ for each $\xi$. Then, solving this ODE system by a modified
variational method, we can show that $\lambda(|\xi|)>0$ is continuous
function on $(0,\infty)$, and the normal modes with spatial frequency
grow in time, providing thus a mechanism for the global linear RT
instability. Consequently, we form a Fourier synthesis of the normal
mode solutions constructed for each spatial frequency $\xi$ to
construct solutions of the linearized equations that
grow in time, when measured in $H^k(\mathbb{R}^3)$ for any $k\geq
0$. This is the content of Section 2. (ii) In Section 3, we show a
uniqueness result of the linearized problem (see Theorem \ref{thm:0301})
in the sense of strong solutions. In spite of the uniqueness, the
linearized problem is global unstable in $H^k(\mathbb{R}^3)$ for any $k$.
(iii) Then we derive some nonlinear energy estimates of the
perturbed problem with small initial data, which make it possible
to take to the limit in the scaled perturbed problem to obtain
the corresponding linearized equations. (iv) Finally, in Section 5,
with the help of the results established in Sections 2--4 and the
Lipschitz structure of $F$, we can obtain the instability of the
 nonlinear problem in the sense of (\ref{0111}). In the proof,
 we shall see that the stability structure (\ref{0114}) would give rise to
certain estimates of solutions to the linearized problem
(\ref{0106})--(\ref{0108}) that cannot hold in general because of
Theorem \ref{thm:0201}.

We should point out that the RT instability based on the Lipschitz
structure was studied by Guo and Tice in \cite{GYTI1} for
compressible inviscid fluids, where the instability was shown in the
$H^3$-norm of $(\varrho ,\mathbf{u})$ and the flow map
 (see (\ref{0512}) in \cite{GYTI1}). Our instability result
Theorem \ref{thm:0101} differs from that of Guo and Tice in that
only $\|u_3\|_{L^2}$ is needed here to describe the instability. This is
also different from that of Guo and Hwang \cite{HHJGY}, in which the
instability for an inhomogeneous incompressible inviscid
fluid is described by the norm $\|(\varrho ,\mathbf{u})\|_{L^2(\Omega)}$.
Roughly speaking, our instability in terms of $\|u_3\|_{L^2}$ only
is based on two important observations: (i) one can construct a
solution $(\varrho,\mathbf{u})$ with $\|u_3(0)\|_{L^2}>0$ of the
linearized problem; (ii) more regularity of the solution $(\varrho,\mathbf{u})$
to the corresponding nonlinear problem can be derived from the problem
(\ref{0105})--(\ref{0107}), we refer to Section 5 for details.
We also mention that in the current paper we have to employ new techniques to construct growing
in time solutions to the linearized problem. To construct such solutions, we shall first transform
the linearized equations to an ODE system. In the the inviscid fluid case, the ODE system can be
viewed as an eigenvalue problem with eigenvalue $-\lambda^2$ (cf.
ODE ({10}) in \cite{HHJGY} for impressible fluids, or ODE system
({3.11}) in \cite{GYTI1} for compressible fluids). Unfortunately,
when the viscosity is present, the linear term multiplied by
$\lambda$ breaks down the natural variational structure, such that
the variational method can not be used. In order to circumvent this
problem, for compressible viscous fluids, recently Guo and Tice
\cite{GYTI2} artificially removed the linear dependence on $\lambda$
by first defining $s:={\lambda}>0$, then solving the family of
modified problems for each $s>0$, and finally showing $s=\lambda(s)$
for some $s$. In \cite{GYTI2} the ODE system is defined in a bounded
domain, and the compact imbedding, an important step in their
construction, can thus be applied. In our case, however, 
our ODE  (see (\ref{0207})) is defined on $\mathbb{R}$, and
consequently the compact embedding does not hold. To overcome this
difficult, here we exploit the property of weak convergence and the
structure of the energy functional $E(\psi)$ corresponding to our
ODE. In particular, we develop a new analysis technique to prove
$\lambda(|\xi|)\in C^0(0,\infty)$ by first showing
$\lambda(|\xi|,s)\in C^0(0,\infty)$ for each fixed $s$, and then
exploiting the monotonicity of $\lambda(\cdot,s)$ to further verify
$\lambda(|\xi|):=\lambda(|\xi|,\lambda)\in C^0(0,\infty)$, see the
proof of Proposition \ref{pro:0205} for details.

We end this section by briefly reviewing some of the
previous results on the nonlinear RT instability for two layer incompressible
fluids with a free interface, where the RT steady state solution is
a denser fluid lying above a lighter one separated by a free interface.
When the densities of two layer fluids are two constants,
Pr$\mathrm{\ddot{u}}$ess and Simonett \cite{JPGS09} used the
$C^0$-semigroup theory and the Henry instability theorem to show the (local)
existence of nonlinear unstable solutions in the
Sobolev-Slobodeckii spaces, where the instability term is described
by the sum of $\|\mathbf{u}\|_{W_p^{2-2/p}}$ and
$\|\mathbf{h}\|_{W_p^{3-2/p}}$ (see \cite[Theorem 1.2]{JPGS09} for
details). When densities of two layer fluids are variable,
 to our best knowledge, the (local) existence of solutions to the
nonlinear problem is still not known unfortunately, and thus the nonlinear
instability is still open. For compressible fluids there are very few results on
the nonlinear RT instability. Guo and Tice proved the instability
of immiscible compressible inviscid fluids in the frame of Lagrangian coordinates
under the assumption of the existence of solutions \cite{GYTI1}, which is
in some sense a compressible analogue to the local
ill-posedness of the RT problem for incompressible fluids given in
\cite{EDGTC1111}.
Recently, Jiang, Jiang and Wang \cite{JFJSWYJ} adapted Guo and
Tice's approach to investigate the nonlinear instability of two
immiscible incompressible fluids with or without surface tension in
Eulerian coordinates without the help of a coordinate
transformation. We remark that the analogue of the RT instability
arises when the fluids are electrically conducting and a magnetic
field is present, and the growth of the instability will be
influenced by the magnetic field due to the generated
electromagnetic induction and the Lorentz force. Some authors have
extended the partial results concerning the RT instability of
superposed flows to the case of MHD flows by overcoming the more
complicated structure due to presence of the magnetic field, see
\cite{HHVQ,KMSMSP,WYC,JFJSWWWO,RDFJSJO}.

\textit{Notation:} Throughout this article we shall repeatedly use the abbreviations:
\begin{eqnarray*} &&
W^{m,p}:=W^{m,p}({\mathbb R}^3),\;\; H^m:=H^m({\mathbb R}^3),\;\; L^p:=L^p({\mathbb R}^3),\\
&& \|\cdot\|_{W^{m,p}}:=\|\cdot\|_{W^{m,p}({\mathbb R}^3)},\;\;
\|\cdot\|_{H^m}:=\|\cdot\|_{H^m({\mathbb R}^3)},\;\;
\|\cdot\|_{L^p}:=\|\cdot\|_{L^p({\mathbb R}^3)},\quad \mbox{etc.}
\end{eqnarray*}

\section{Construction of solutions to the linearized problem}
We wish to construct a solution to the linearized equations
(\ref{0108}) that has growing $H^k$-norm for any $k$. We will
construct such solutions via Fourier synthesis by first constructing
a growing mode for any but fixed spatial frequency.

\subsection{Linear growing modes}

To begin, we make a growing mode ansatz of solutions, i.e.,
\begin{equation*}
{\varrho}(\mathbf{x})=\tilde{\rho} (\mathbf{x})e^{\lambda t},\;\;
\mathbf{u}(\mathbf{x})=\tilde{\mathbf{v}}(\mathbf{x})e^{\lambda t},\;\;
{q}(\mathbf{x})=\tilde{p}(\mathbf{x})e^{\lambda t}\quad\mbox{for
some }\lambda>0.
\end{equation*}
Substituting this ansatz into (\ref{0108}), and then eliminating
$\tilde{\rho}$ by using the first equation, we arrive at the
time-invariant system for
$\tilde{\mathbf{v}}=(\tilde{v}_1,\tilde{v}_2,\tilde{v}_3)$ and
$\tilde{p}$:
\begin{equation}\label{0201}
\left\{
                              \begin{array}{ll}
\lambda^2\bar{\rho}\tilde{\mathbf{v}}+\lambda\nabla \tilde{{p}}
=\lambda\mu \Delta \tilde{\mathbf{v}}+g\bar{\rho}'\tilde{v}_3e_3,\\[1mm]
\mathrm{div}\, \tilde{\mathbf{v}}=0
\end{array}
                            \right.
\end{equation}
with
\begin{equation*}\label{0202}
\lim_{|\mathbf{x}|\rightarrow
+\infty}\tilde{\mathbf{v}}(\mathbf{x})=\mathbf{0}.
\end{equation*}

 We fix a spatial frequency $\xi=(\xi_1,\xi_2)\in \mathbb{R}^2$
and take the horizontal Fourier transform
of $(\tilde{v}_1,\tilde{v}_2,\tilde{v}_3)$ in (\ref{0201}), which we
denote with either $\hat{\cdot}$ or $\mathcal{F}$, i.e.,
\begin{equation*}\hat{f}(\xi,x_3)=\int_{\mathbb{R}^2}f(x',x_3)e^{-ix'\cdot
\xi}\mathrm{d}x'.\end{equation*} Define the new unknowns
\begin{equation*}\label{0203}\varphi(x_3)=i\hat{\tilde{v}}_1(\xi,x_3), \;\;
\theta(x_3)=i\hat{\tilde{v}}_2(\xi,x_3),\;\;\psi(x_3)=\hat{\tilde{v}}_3(\xi,x_3),\;\;
\pi(x_3)={\hat{\tilde{p}}}(\xi,x_3),\end{equation*} so that
\begin{equation*}\label{0204}
\mathcal{F}(\mathrm{div}\tilde{\mathbf{v}})=\xi_1\varphi+\xi_2\theta+\psi',
\end{equation*}
where $'=d/dx_3$. Then, for $\varphi$, $\theta$, $\psi$ and
$\lambda=\lambda(\xi)$ we arrive at the following system of ODEs.
\begin{equation}\label{0205} \left\{
                              \begin{array}{ll}
\lambda^2\bar{\rho} \varphi-\lambda\xi_1\pi+\lambda\mu (|\xi|^2\varphi-\varphi'')=0,\\[1mm]
\lambda^2\bar{\rho} \theta-\lambda\xi_2\pi+\lambda\mu (|\xi|^2\theta-\theta'')=0,\\[1mm]
\lambda^2\bar{\rho} \psi+\lambda\pi'+\lambda\mu (|\xi|^2\psi-\psi'')=g\bar{\rho}'\psi,\\[1mm]
\xi_1\varphi+\xi_2\theta+\psi'=0
\end{array}
                            \right.
\end{equation}
 with
\begin{eqnarray}\label{0206}
\varphi(-\infty)=\theta(-\infty)=\psi(-\infty)=\varphi(+\infty)=\theta(+\infty)=\psi(+\infty)=0.
\end{eqnarray}

Eliminating $\pi$ from the third equation in (\ref{0205}) we obtain
the following ODE for $\psi$
\begin{equation}\label{0207}
-\lambda^2[\bar{\rho}|\xi|^2\psi-(\bar{\rho}\psi')']=\lambda\mu
(|\xi|^4\psi-2|\xi|^2\psi''+\psi'''')-g|\xi|^2\bar{\rho}'\psi
\end{equation}
with
\begin{eqnarray}\label{0208}
&&\psi(-\infty)=\psi'(-\infty)=\psi(+\infty)=\psi'(+\infty)=0.
\end{eqnarray}

Similarly to \cite{GYTI2}, we can apply the variational method to
construct a solutions of (\ref{0207})--(\ref{0208}). The idea of the
proof can be found in Guo and Tice's paper for viscous compressible
flows \cite{GYTI2}, and was adapted later by other authors to
investigate the instability for other fluid models
\cite{JFJSWWWO,JJHTI,WYC,RDFJSJS}.

Now we fix a non-zero vector $\xi\in \mathbb{R}^2$ and $s>0$. From
(\ref{0207})--(\ref{0208}) we get a family of the modified problems
\begin{equation}\label{0209}
-\lambda^2[\bar{\rho}|\xi|^2\psi-(\bar{\rho}\psi')']=s\mu
(|\xi|^4\psi-2|\xi|^2\psi''+\psi'''')-g|\xi|^2\bar{\rho}'\psi,
\end{equation}
coupled with (\ref{0208}). We define the energy functional of
(\ref{0209}) by
\begin{equation}\label{0210}E(\psi)=\int_{\mathbb{R}}s\mu(4|\xi|^2|\psi'|^2
+||\xi|^2\psi+\psi''|^2)-g|\xi|^2
\bar{\rho}'\psi^2\mathrm{d}x_3\end{equation} with a associated
admissible set
\begin{equation}\label{0211}
\mathcal{A}=\left\{\psi\in
H^2(\mathbb{R})~\bigg|~J(\psi):=\int_{\mathbb{R}}\bar{\rho}(|\xi|^2|\psi|^2+|\psi'|^2)\mathrm{d}x_3=1\right\}.
\end{equation}
 Thus we can find a $-\lambda^2$ by minimizing
\begin{equation}\label{0212}
-\lambda^2(|\xi|)=\alpha(|\xi|):=\inf_{\psi\in
\mathcal{A}}E(\psi).\end{equation}
In order to emphasize the dependence on $s\in (0,\infty)$ we will sometimes write
\begin{equation*}
E(\psi,s):=E(\psi)\mbox{ and } \alpha(s):=\inf_{\psi\in
\mathcal{A}}E(\psi,s).\end{equation*}

Next we show that a minimizer
of (\ref{0212}) exists for the case of
$\inf_{\mathcal{A}}E(\psi,s)<0$, and that the corresponding 
Euler-Lagrange equations are equivalent to (\ref{0208}), (\ref{0209}).
\begin{pro}\label{pro:0201}
 For any fixed $\xi$ with $|\xi|\neq 0$, $\inf_{\psi\in
\mathcal{A}}E(\psi,s)>-\infty$. In particular, if there exists a
$\bar{\psi}\in \mathcal{A}$, such that $E(\bar{\psi})<0$, then $E$
achieves its infinimum on $\mathcal{A}$. In addition, let $\psi$ be
a minimizer and $-\lambda^2:=E(\psi)$, then the pair ($\psi$,
$\lambda^2$) satisfies (\ref{0208}), (\ref{0209}). Moreover,
$\psi\in H^k(\mathbb{R})$ for any positive integer $k$.
\end{pro}
\begin{pf}
 We first note that for any $\psi\in \mathcal{A}$,
\begin{equation}\label{0213}E(\psi)\geq -{g|\xi|^2}\int_{\mathbb{R}}
\bar{\rho}'\psi^2\mathrm{d}x_3\geq
-{g}\left\|\frac{\bar{\rho}'}{\bar{\rho}}\right\|_{L^\infty(\mathbb{R})}\int_{\mathbb{R}}
\bar{\rho}|\xi|^2\psi^2\mathrm{d}x_3\geq
-{g}\left\|\frac{\bar{\rho}'}{\bar{\rho}}\right\|_{L^\infty(\mathbb{R})}.
\end{equation}
Hence $E$ is bounded from below on $\mathcal{A}$ by virtue of (\ref{0102}). Let $\psi_n\in
\mathcal{A}$ be a minimizing sequence, then $E(\psi _n)$ is bounded.
This together with (\ref{0211}) and (\ref{0210}) again implies that
$\psi_n$ is bounded in $H^2(\mathbb{R})$. So, there exists a $\psi\in
H^2(\mathbb{R})$, such that $\psi_n\rightarrow \psi$ weakly in $H^2(\mathbb{R})$
and strongly in $H^1_{\mathrm{loc}}(\mathbb{R})$.
Moreover, by the lower semi-continuity, locally strong convergence,
  (\ref{0102}) and the assumption that $E(\bar{\psi})<0$ for some
$\bar{\psi}\in \mathcal{A}$, we have
\begin{equation*}E(\psi)\leq \liminf_{n\rightarrow
\infty}E(\psi_n)=\inf_{\mathcal{A}}E<0,\quad\mbox{and }\; 0<J(\psi)\leq 1.
\end{equation*}

Suppose by contradiction that $J(\psi)<1$. By the homogeneity
of $J$ we may find an $\alpha>1$ so that $J(\alpha \psi)=1$, i.e., we
may scale up $\psi$ so that $\alpha \psi\in{\mathcal{A}}$.  From
this we deduce that
\begin{eqnarray*}
E(\alpha \psi)=\alpha^2E(\psi)\leq
\alpha^2\inf\limits_{\mathcal{A}}E<\inf\limits_{\mathcal{A}}E<0,
\end{eqnarray*}
which is a contradiction since $\alpha \psi\in{\mathcal{A}}$. Hence
$J(\psi)=1$ so that $\psi\in{ \mathcal{A}}$. This shows that $E$
achieves its infinimum on $\mathcal{A}$.

Notice that since $E$ and $J$ are homogeneous of degree 2,
(\ref{0212}) is equivalent to
  \begin{equation}\label{0215}\alpha(s)=\inf_{\psi\in
  H^2(\mathbb{R})}\frac{E(\psi)}{J(\psi)}.
\end{equation}
For any $\tau\in \mathbb{R}$ and $\psi_0\in H^2(\mathbb{R})$
we take $\psi(\tau)=\psi+\tau\psi_0$, then (\ref{0215}) implies
  \begin{equation*}E(\psi(\tau))+\lambda^2J(\psi(\tau))\geq 0.
\end{equation*}
If we set $I(\tau)=E(\psi(\tau))+\lambda^2J(\psi(\tau))$, then we see that
$I(\tau)\geq 0$ for all $\tau\in \mathbb{R}$ and $I(0)=0$. This
implies $I'(0)=0$. By virtue of (\ref{0210}) and (\ref{0211}), a
direct computation leads to
  \begin{equation}\label{0217}\begin{aligned}
  & s\mu \int_{\mathbb{R}}(4|\xi|^2\psi'\psi'_0+(|\xi|^2\psi+\psi'')
  (|\xi|^2\psi_0+\psi_0''))\mathrm{d}x_3\\
 & =   g|\xi|^2\int_\mathbb{R}\bar{\rho}'\psi\psi_0\mathrm{d}x_3
  -\lambda^2\int_{\mathbb{R}}\bar{\rho}(|\xi|^2\psi\psi_0+\psi'\psi_0')\mathrm{d}x_3,
\end{aligned}\end{equation}
where we have used the upper boundedness of $\bar{\rho}$.

 By further assuming that $\psi_0$ is compactly supported in $\mathbb{R}$, we
find that $\psi$ satisfies the equation (\ref{0209}) in the weak sense
on $\mathbb{R}$. In order to improve the regularity of $\psi$, we rewrite (\ref{0217}) as
  \begin{equation}\label{0218}\begin{aligned}
 \int_\mathbb{R}
\psi''\psi_0''\mathrm{d}x_3=&\frac{1}{s\mu}\int_{\mathbb{R}}\left(g|\xi|^2\bar{\rho}'\psi-
\lambda^2(|\xi|^2\bar{\rho}\psi-(\bar{\rho}\psi')')+s\mu(2|\xi|^2\psi''
-|\xi|^4\psi) \right)\psi_0\mathrm{d}x_3\\
=&\int_{\mathbb{R}}f \psi_0\mathrm{d}x_3.
\end{aligned}\end{equation}

For any $n\geq 1$, let $\psi_{1,n},\psi_2\in C_0^\infty(\mathbb{R})$
satisfy $\psi_{1,n}(x_3)\equiv 1$ for $|x_3|\leq n$. If we take
$\psi_0=\psi_{1,n}\int_{-\infty}^{x_3}\psi_2\mathrm{d}\tau$ in (\ref{0218}), then we have
  \begin{equation*}\begin{aligned}\int_\mathbb{R}
(\psi_{1,n}\psi'')\psi_2'\mathrm{d}x_3=&
\int_{\mathbb{R}}\left(f\psi_{1,n}\int_{-\infty}^{x_3}\psi_2\mathrm{d}\tau
-\psi_{1,n}''\psi''\int_{-\infty}^{x_3}\psi_2\mathrm{d}\tau-2\psi'_{1,n}\psi''
\psi_2\right)\mathrm{d}x_3 \\
= & \int_{\mathbb{R}}\left(\int_{x_3}^{+\infty}(f\psi_{1,n}-\psi_{1,n}''\psi'')\mathrm{d}\tau-2\psi'_{1,n}\psi''
\right)\psi_2\mathrm{d}x_3,
\end{aligned}\end{equation*}
which, recalling $\psi\in H^2(\mathbb{R})$, implies $\psi''\in H^1_{\mathrm{loc}}(\mathbb{R})$ and
\begin{equation*}\begin{aligned}
\psi'''=(\psi_{1,n}\psi'')'=\int_{x_3}^{+\infty}(f\psi_{1,n}-\psi_{1,n}''\psi'')\mathrm{d}\tau\quad
 \mbox{ for any }x_3 \mbox{ with }|x_3|\leq n.
\end{aligned}\end{equation*}
Integrating by parts, we can rewrite (\ref{0218}) as
  \begin{equation*}\begin{aligned}
& -\int_\mathbb{R}
\psi'''\psi_0'\mathrm{d}x_3=\frac{1}{s\mu}\int_{\mathbb{R}}\left(g|\xi|^2\bar{\rho}'\psi-
\lambda^2(|\xi|^2\bar{\rho}\psi-(\bar{\rho}\psi')')+s\mu(2|\xi|^2\psi''
-|\xi|^4\psi) \right)\psi_0\mathrm{d}x_3,
\end{aligned}\end{equation*}
which, keeping in mind that $\psi\in H^2(\mathbb{R})$, yields
$\psi''''\in L^2(\mathbb{R})$. Hence $\psi\in H^4_{\mathrm{loc}}(\mathbb{R})\cap
C^{3,1/2}_{\mathrm{loc}}(\mathbb{R})$, and
$\psi^{'}(\infty)=\psi^{''}(\infty)=\psi^{'''}(\infty)=0$. Using
these facts, H\"{o}lder's inequality, and integration by
parts, we deduce that
  \begin{equation}\begin{aligned}\label{0221}
\|\psi'''\|_{L^2(\mathbb{R})}^2=\int_{\mathbb{R}}|\psi'''|^2\mathrm{d}x_3
=-\int_{\mathbb{R}}\psi''\psi''''\mathrm{d}x_3\leq\|\psi''\|_{L^2(\mathbb{R})}\|\psi''''\|_{L^2(\mathbb{R})},
\end{aligned}\end{equation}
i.e., $\psi'''\in L^2(\mathbb{R})$. Consequently, $\psi\in H^4(\mathbb{R})$
and solves (\ref{0208})--(\ref{0209}). This immediately gives
that $\psi\in H^k(\mathbb{R})$ for any positive integer $k\geq 5$. \hfill $\Box$
\end{pf}

Next, we want to show that there is a fixed point such that
$\lambda=s$. To this end, we first give some properties of
$\alpha(s)$ as a function of $s> 0$.

\begin{pro}\label{pro:0202} The function $\alpha(s)$ defined on $(0,\infty)$ enjoys the following
properties:
\begin{enumerate}[\quad \ (1)]
 \item  For any $a$, $b\in (0,\infty)$ with $a<b$, there exist constants $c_1$, $c_2>0$ depending on $\bar{\rho}$,
  $\mu$, $g$, $a$ and $b$, such that
  \begin{equation}\label{0222}\alpha(s)\leq -c_1+sc_2\quad\mbox{ for all }|\xi|\in [a, b].\end{equation}
  \item $\alpha(s)\in C_{\mathrm{loc}}^{0,1}(0,\infty)$ is nondecreasing.
\end{enumerate}
\end{pro}
\begin{pf} (1)  In view of (\ref{0103}), there exists a $\tilde{\psi}\in
C_0^\infty(\mathbb{R})$ such that
  \begin{equation}\label{0224} \frac{
  ga^2\int_{\mathbb{R}}
  \bar{\rho}'\tilde{\psi}^2\mathrm{d}x_3}{\int_{\mathbb{R}}\bar{\rho}(b^2|\tilde{\psi}|^2+|\tilde{\psi}'|^2)\mathrm{d}x_3}:=c_1>
  0,
\end{equation}
where the constant $c_{1}$ depends on $a$, $b$, $g$ and
$\bar{\rho}$. Now, we use (\ref{0215}) and (\ref{0224}) to find that
  \begin{equation*}\begin{aligned}\alpha(s)&=\inf_{\psi\in
  H^2(\mathbb{R})}\frac{E(\psi)}{J(\psi)}\leq \frac{\int_{\mathbb{R}}s
  \mu(4|\xi|^2|\tilde{\psi}'|^2 +||\xi|^2\tilde{\psi}+\tilde{\psi}''|^2)-
g|\xi|^2\bar{\rho}'\tilde{\psi}^2\mathrm{d}x_3}{\int_{\mathbb{R}}
\bar{\rho}(|\xi|^2|\tilde{\psi}|^2+|\tilde{\psi}'|^2)\mathrm{d}x_3} \\[1mm]
&\leq  s \frac{\mu\int_{\mathbb{R}}(4|\xi|^2|\tilde{\psi}'|^2
+||\xi|^2\tilde{\psi}+\tilde{\psi}''|^2)\mathrm{d}x_3}{\int_{\mathbb{R}}
\bar{\rho}(|\xi|^2|\tilde{\psi}|^2+|\tilde{\psi}'|^2)\mathrm{d}x_3}- \frac{ga^2\int_{\mathbb{R}}
  \bar{\rho}'\tilde{\psi}^2\mathrm{d}x_3}{\int_{\mathbb{R}}\bar{\rho}
  (b^2|\tilde{\psi}|^2+|\tilde{\psi}'|^2)\mathrm{d}x_3} := sc_2-c_1,
\end{aligned}\end{equation*}
where the positive constant $c_2$ depends on $\bar{\rho}$,
  $\mu$, $g$, $a$ and $b$. Hence, (\ref{0222}) holds.

(2) To show the second assertion, we let $Q:=[a,b]\subset
\mathbb{R}$ be a bounded interval, and
$$E_1(\psi)=\mu\int_{\mathbb{R}}(4|\xi|^2|\psi'|^2 +|\,|\xi|^2\psi+\psi''|^2)\mathrm{d}x_3.$$
For any $s\in  Q $, there exists a minimizing sequence $\{\psi^n_{s}\}\subset \mathcal{A}$ of
$\inf_{\psi\in \mathcal{A}}E(\psi,s)$, such that
\begin{equation}\label{0225}\begin{aligned}|\alpha(s)-E(\psi_{s}^n,s)|<1
\end{aligned}.\end{equation}
  Making use of (\ref{0210}), (\ref{0213}),
(\ref{0222}) and  (\ref{0225}), we infer that
\begin{equation}\begin{aligned}\label{0226}E_1(\psi_s^n,s)=&\frac{E(\psi_s^n,s)}{s}
+ \frac{g|\xi|^2}{s}\int_{\mathbb{R}} \bar{\rho}'\psi^2\mathrm{d}x_3  \\
\leq & \frac{1+\max\{|bc_2-c_1|,{g}\left\|{\bar{\rho}'}/{\bar{\rho}}
\right\|_{L^\infty(\mathbb{R})}\}}{a}+\frac{g}{a}\left\|\frac{\bar{\rho}'}{\bar{\rho}}
\right\|_{L^\infty(\mathbb{R})}:=K.
\end{aligned}\end{equation}

For $s_i\in Q$ ($i=1,2$), we find that
\begin{equation}\begin{aligned}\label{0227}\alpha(s_1)\leq \limsup_{n\rightarrow
\infty}E(\varphi_{s_2}^n,s_1)\leq & \limsup_{n\rightarrow
\infty}E(\psi_{s_2}^n,s_2)+|s_1-s_2|\limsup_{n\rightarrow
\infty}E_1(\psi_{s_2}^n)\\
\leq & \alpha(s_2)+K|s_1-s_2|,
\end{aligned}\end{equation}where  $\{\psi^n_{s_2}\}\subset\mathcal{A}$ is a minimizing sequence of
$\inf_{\psi\in \mathcal{A}}E(\psi,s_2)$ and the constant $K$ is given in (\ref{0226}).
Reversing the role of the indices 1 and 2 in the
derivation of the inequality (\ref{0227}), we obtain the same boundedness with
the indices switched. Therefore, we deduce that
\begin{equation*}\begin{aligned}|\alpha(s_1)-\alpha(s_2)|\leq K|s_1-s_2|,
\end{aligned}\end{equation*}
which yields $\alpha(s)\in C_{\mathrm{loc}}^{0,1}(0,\infty)$.

Finally, from (\ref{0210}) and (\ref{0212}) it follows that
\begin{equation*}
\alpha(s_1)\leq \limsup_{n\rightarrow\infty}E(\psi_{s_2}^n,s_1)
\leq \limsup_{n\rightarrow\infty}E(\psi_{s_2}^n,s_2)=\alpha(s_2)\;
\mbox{ for any }0<s_1<s_2<\infty.
\end{equation*}
Hence $\alpha(s)$ is nondecreasing on $(0,\infty)$.
This completes the proof of Proposition \ref{pro:0202}.
 \hfill $\Box$
\end{pf}

Given $\xi\in \mathbb{R}^2$ with $|\xi|\neq 0$, by virtue of
(\ref{0222}), there exists a $s_0>0$ depending on the quantities $\bar{\rho}$,
$\mu$, $g$, $|\xi|$, such that for any $s\leq s_0$, $\alpha(s)<0$. Let
\begin{equation}\label{0230}\mathfrak{S}_{|\xi|}:=\sup\{s~|~\alpha(\tau)<0\mbox{ for any }\tau\in
(0,s)\}>0,\end{equation}  then $\mathfrak{S}_{|\xi|}>0$. This allows
us to define $\lambda(s)=\sqrt{-\alpha(s)}>0$ for any $s\in
\mathcal{S}_{|\xi|}:=(0,\mathfrak{S}_{|\xi|})$. Therefore, as a
result of Proposition \ref{pro:0201}, we have the following existence for the modified problem
(\ref{0209}), (\ref{0210}).
\begin{pro}\label{pro:0203}
For each $|\xi|\neq 0$ and $s\in \mathcal{S}_{|\xi|}$ there is
a solution $\psi=\psi(|\xi|,x_3)\equiv\!\!\!\!\!\!/\ 0$ with
$\lambda=\lambda(|\xi|,s)>0$ to the problem (\ref{0208}), (\ref{0209}).
Moreover, $\psi\in H^k(\mathbb{R})$ for any positive integer $k$.
\end{pro}

 Now, we can use Proposition \ref{pro:0202}, (\ref{0230}) and (\ref{0213})
to check that $\lambda(s)\in C_{\mathrm{loc}}^{0,1}(\mathcal{S}_{|\xi|})$ is nonincreasing
(in fact, we can further show that $\lambda(s)$ is strictly increasing, we refer to the
proof of \cite[Proposition 3.6]{GYTI2}),
$\lambda(s)\leq \sqrt{g}\|\sqrt{{\bar{\rho}'}/{\bar{\rho}}}\|_{L^\infty(\mathbb{R})}$,
and $\lim_{|\xi|\rightarrow \mathfrak{S}_{|\xi|}}\lambda(\mathfrak{S}_{|\xi|})=0$
if $\mathfrak{S}_{|\xi|}<+\infty$.
Hence, we can employ a fixed-point argument to find $s\in\mathcal{S}_{|\xi|}$ so that
$s=\lambda(|\xi|, s)$, thus and obtain a
solution to the original problem (\ref{0207}), (\ref{0208}).
\begin{pro}\label{pro:0204} Let $|\xi|\neq 0$,
then there exists a unique $s\in \mathcal{S}_{|\xi|}$, such that
$\lambda(|\xi|,s)=\sqrt{-\alpha(s)}>0$ and $s=\lambda(|\xi|,s)$.
\end{pro}
\begin{pf}
We refer to \cite[Theorem 3.8]{GYTI2} (or \cite[Lemma 3.7]{WYC}) for
a proof.\hfill $\Box$
\end{pf}

 Consequently, in view of Propositions \ref{pro:0203} and \ref{pro:0204},
 we conclude the following existence for the problem (\ref{0207}), (\ref{0208}).
\begin{thm}\label{thm:0201}
For each $|\xi|\neq 0$, there exist
$\psi=\psi(|\xi|,x_3)\equiv\!\!\!\!\!\!/\ 0$ and $\lambda(|\xi|)>0$
satisfying (\ref{0207}), (\ref{0208}). Moreover, $\psi\in H^k(\mathbb{R})$ for
any positive integer $k$.
\end{thm}

We end this subsection by giving some properties of the solutions
established in Theorem \ref{thm:0201} in terms of $\lambda(|\xi|)$,
which show that $\lambda$ is a bounded, continuous function of $|\xi|$.
\begin{pro}\label{pro:0205}
The function
$\lambda :(0,\infty)\rightarrow (0,\infty)$ is continuous and satisfies
\begin{equation}\label{0231}
\sup_{0<|\xi|<\infty}\lambda(|\xi|)\leq
\sqrt{g}\left\|\sqrt{{\bar{\rho}'}/{\bar{\rho}}}\right\|_{L^\infty(\mathbb{R})}.\end{equation}
\end{pro}
\begin{pf} The boundedness of $\lambda$ (\ref{0231}) follows
from (\ref{0213}).  To show the continuity of $\lambda$, we see that for any but fixed
$\xi_0\not =0$, there exists an interval $[a,b]\subset (0,\infty)$ so that $|\xi_0|\in (a,b)$. Assume
 $|\xi|\rightarrow |\xi_0|$ with $|\xi|\in(a,b)$, and denote $\kappa=|\xi|^2-|\xi_0|^2$, then
$\kappa\rightarrow 0 $ as $|\xi|\rightarrow |\xi_0|$.

(i) We first show
\begin{equation}\label{0234}
\lim_{|\xi|\rightarrow |\xi_0|}
\alpha(|\xi|,s)=\alpha(|\xi_0|,s)\;\;\mbox{ for any }s\in
\mathcal{S}_{|\xi|}.\end{equation}

By virtue of Proposition \ref{pro:0201}, for any $|\xi|\in(a,b)$,
there exists a functions $\psi_{|\xi|}\in \mathcal{A}$, such that
\begin{equation}\label{0235}\alpha(|\xi|)=\int_{\mathbb{R}}s\mu(4|\xi|^2|\psi'_{|\xi|}|^2
+||\xi|^2\psi_{|\xi|}+\psi''_{|\xi|}|^2)-
g\bar{\rho}'|\xi|^2\psi^2_{|\xi|}\mathrm{d}x_3\end{equation}
Utilizing (\ref{0211}) and (\ref{0222}), we have
\begin{equation}\label{0236}\|\psi_{|\xi|}\|_{H^2(\mathbb{R})}\leq  c_5,\end{equation}
where $c_5$ depends on $\bar{\rho}$,  $\mu$, $g$, $a$, $b$ and $s$.

 Substitution of $|\xi|^2=|\xi|^2+\kappa$ into (\ref{0235}) results in
\begin{equation}\label{0237}\begin{aligned}\alpha(|\xi|)=&\int_{\mathbb{R}}s\mu(4|\xi_0|^2|\psi'_{|\xi|}|^2
+||\xi_0|^2\psi_{|\xi|}+\psi''_{|\xi|}|^2)-
g|\xi_0|^2\bar{\rho}'\psi^2_{|\xi|}\mathrm{d}x_3+\kappa
f(\kappa,\psi_{|\xi|})\\
\geq &\alpha(|\xi_0|)+\kappa
f(\kappa,\psi_{|\xi|}),\end{aligned}\end{equation} where
\begin{equation*}\label{0238}f(\kappa,\psi_{|\xi|})=\int_{\mathbb{R}}s\mu(4|\psi'_{|\xi|}|^2 +2
\psi_{|\xi|}(|\xi|^2\psi_{|\xi|}+\psi''_{|\xi|})+\kappa\psi^2_{|\xi|}-
g|\xi|^2\bar{\rho}'\psi^2_{|\xi|}\mathrm{d}x_3.\end{equation*}
By H\"{o}lder's inequality and (\ref{0236}), we can bound
\begin{equation}\label{0239}\begin{aligned}|
f(\kappa,\psi_{|\xi|})|\leq c_6\mbox{ for some constant
}c_6.\end{aligned}\end{equation}

Similarly to (\ref{0237}) and (\ref{0239}), we also have
\begin{equation}\label{0240}\begin{aligned}\alpha(|\xi_0|)\geq\alpha(|\xi|)-\kappa
f(-\kappa,\psi_{|\xi_0|})\mbox{ and }|f(-\kappa,\psi_{|\xi_0|})|\leq
c_6.\end{aligned}\end{equation} Combining (\ref{0237}) with
(\ref{0240}), we  get
\begin{equation*}\label{0241}\begin{aligned}\kappa
f(-\kappa,\psi_{|\xi_0|})\geq \alpha(|\xi|)-\alpha(|\xi_0|)\geq
\kappa f(\kappa,\psi_{|\xi|}),\end{aligned}\end{equation*}
which, together with (\ref{0239}) and (\ref{0240}), implies that (\ref{0234}). Hence
\begin{equation}\label{0242}
\lim_{|\xi|\rightarrow |\xi_0|}
\lambda(|\xi|,s)=\lambda(|\xi_0|,s)\mbox{ for any }s\in
\mathcal{S}_{|\xi|},\end{equation}
because of $\lambda(|\xi|,s)=\sqrt{-\alpha(|\xi|,s)}$.

(ii) Exploiting (\ref{0242}) and Propositions \ref{pro:0204}, we
know that for any $\varepsilon>0$, there exists a $\delta>0$ such that
 $|\lambda(|\xi|,s_{|\xi_0|})-\lambda(|\xi_0|,s_{|\xi_0|})|<\varepsilon$  and
 $s_{|\xi_0|}=\lambda({|\xi_0|},s_{|\xi_0|})=\sqrt{-\alpha({|\xi_0|,s_{|\xi_0|}})}$ for any
 $|\,|\xi|-|\xi_0|\,|<\delta$.
 On the other hand, for each $|\xi|>0 $, $\lambda(s)$ is nonincreasing and continuous
 on $\mathcal{S}_{|\xi|}$, and there exists a unique
 $s_{|\xi|}\in \mathcal{S}_{|\xi|}$ satisfying
$\lambda(|\xi|,s_{|\xi|})=s_{|\xi|}>0$ by Propositions
\ref{pro:0204}. Consequently, we immediately infer that
$|\lambda(|\xi|,s_{|\xi|})-\lambda(|\xi_0|,s_{|\xi_0|})|<\varepsilon$
with  $s_{|\xi|}=\lambda(|\xi|,s_{|\xi|})$. Hence $\lambda(|\xi|)$
is continuous.  This completes the proof of the proposition.
 \hfill $\Box$
\end{pf}
\begin{rem}In addition, since $\bar{\rho}'\in C_0^\infty(\mathbb{R})$, we can bound
$\lambda$ as follows.

(i) Applying integrating by parts and H$\mathrm{\ddot{o}}$lder
inequality,
\begin{equation*}\begin{aligned} \lambda^2(|\xi|)\leq
&\int_{\mathbb{R}} g|\xi|^2\bar{\rho}'\psi^2_{|\xi|}\mathrm{d}x_3=2
{g|\xi|}\int_{\mathbb{R}}
\bar{\rho}|\xi|\psi_{|\xi|}\psi_{|\xi|}'\mathrm{d}x_3\\
\leq &2{g|\xi|}\left(\int_{\mathbb{R}}
\bar{\rho}|\xi|^2\psi_{|\xi|}^2\mathrm{d}x_3\right)^{\frac{1}{2}}\left(\int_{\mathbb{R}}
\bar{\rho}|\psi_{|\xi|}'|^2\mathrm{d}x_3\right)^{\frac{1}{2}}\leq
2{g|\xi|}.
\end{aligned}\end{equation*}
Consequently we have $\lim_{|\xi|\rightarrow 0}\lambda(|\xi|)=0$.

 (ii) There exists a functions $\psi_{|\xi|}\in \mathcal{A}$ such that
\begin{equation*}-\lambda^2(|\xi|)=\lambda(|\xi|)\mu\int_{\mathbb{R}}(4|\xi|^2|\psi'_{|\xi|}|^2
+||\xi|^2\psi_{|\xi|}+\psi''_{|\xi|}|^2)-
g\bar{\rho}'|\xi|^2\psi^2_{|\xi|}\mathrm{d}x_3,\end{equation*} which
implies  that there exists a constant $C({\bar{\rho}})$ depending on
${\bar{\rho}}$ such that
\begin{equation*}   0<4\mu\lambda(|\xi|)\int_{\mathbb{R}}|\psi'_{|\xi|}|^2
\mathrm{d}x_3\leq
gC({\bar{\rho}})\int_{\mathbb{R}}|\psi_{|\xi|}'|^2\mathrm{d}x_3.\end{equation*}
We immediately get that\begin{equation}
\label{n0229n}\lambda(|\xi|)\leq \frac{gC({\bar{\rho}})}{4\mu}.
\end{equation}
\end{rem}

\subsection{Construction of a solution to the system (\ref{0205}), (\ref{0206})}
A solution to (\ref{0207}), (\ref{0208}) gives rise to a solution of
the system  (\ref{0205}), (\ref{0206}) for the growing mode velocity
$\mathbf{u}$ as well.
\begin{thm}\label{thm:0302}
For each $\xi\in \mathbb{R}^2$ with $|\xi|>0$, there exists a
solution $({\varphi},{\theta},{\psi},{\pi})= ({\varphi}(\xi,
x_3),{\theta}(\xi, x_3),$ ${\psi}(|\xi|,x_3), {\pi}(|\xi| ,x_3))$
with $\lambda =\lambda(|\xi|)>0$ to (\ref{0205}), (\ref{0206}), and
the solution belongs to $(H^k(\mathbb{R}))^4$  for any positive integer $k$.
\end{thm}
\begin{pf} With the help of Theorem \ref{thm:0201}, we first construct a solution
$(\psi ,\lambda )=(\psi (|\xi|,x_3),\lambda (|\xi|))$ satisfying
(\ref{0207}), (\ref{0208}). Recalling $\lambda >0$ and
$\psi\in\mathcal{A}\cap H^k(\mathbb{R})$ for any positive integer $k$,
multiplying (\ref{0205})$_1$ and (\ref{0205})$_2$ by $\xi_1$
and $\xi_2$ respectively, adding the resulting equations, and
utilizing (\ref{0205})$_4$, we find that $\pi$ can be expressed by $\psi$, i.e.,
\begin{equation}\label{0243}\pi =\pi(|\xi|,x_3)
=[{\mu \psi'''-(\lambda\bar{\rho}+\mu |\xi|^2)\psi'}]{|\xi|^{-2}}.
\end{equation}

Next, we construct the solution $(\varphi,\theta)$. To this end, we shall exploit the fact that the problem
(\ref{0205}), (\ref{0206}) is invariant under simultaneous rotations
of $(\varphi,\theta)$ and $(\xi_1,\xi_2)$. Indeed, it is easy to see
that if $\mathcal{R}\in SO(2)$ is a rotation operator,
then $\mathcal{R}(\varphi,\theta)$, $\mathcal{R}(\xi_1,\xi_2)$ is
also a solution with the same $\psi$, $\pi$ and $\lambda$. Thus,
given any $\xi$ we choose a rotation operator $\mathcal{R}_{\xi}$ so
that $\mathcal{R}_{\xi}\xi=(|\xi|,0)$. Hence,
\begin{equation*}
(\tilde{\varphi},\tilde{\theta},\psi,\pi,\lambda)=\big(-\psi'_{|\xi|}(|\xi|)/|\xi|,0,\psi(|\xi|),\pi
\mbox{ (given by (\ref{0243}))},\lambda(|\xi|)\big)\end{equation*}
is a solution to (\ref{0205}), (\ref{0206}) with $(\xi_1,\xi_2)=(|\xi|,0)$. Now, if we define
\begin{equation}\label{0244}(\varphi,\theta,\psi,\pi,\lambda):=
(\mathcal{R}_{\xi}^{-1}(\tilde{\varphi},\tilde{\theta}),\psi,\pi,\lambda)=
(-(\xi_1,\xi_2)\psi'_{|\xi|}/|\xi|^2,\psi,\pi,\lambda),\end{equation}
we find that $(\psi,\varphi,\theta,\pi,\lambda)$ constructed
above is indeed a solution to the problem (\ref{0205}), (\ref{0206}).
\hfill $\Box$
\end{pf}
\begin{rem}\label{rem:0201}
For each $x_3$, it is easy to see that the solution
$({\varphi}(\xi,\cdot),{\theta}(\xi,\cdot),{\psi}(|\xi|,\cdot),{\pi}(|\xi|,\cdot),
\lambda(|\xi|))$ constructed in Theorem \ref{thm:0302} has the
following properties:
\begin{enumerate}[\quad \ (1)]
  \item $\lambda(|\xi|)$, ${\psi}(|\xi|,\cdot)$ and ${\pi}(|\xi|,\cdot)$
  are even on $\xi_1$ or $\xi_2$, when the another variable is fixed;
  \item ${\varphi}(\xi,\cdot)$ is odd on $\xi_1$, but even on $\xi_2$,
  when the another variable is fixed;
  \item ${\theta}(\xi,\cdot)$ is even on $\xi_1$, but odd on $\xi_2$,
  when the another variable is fixed.
  \end{enumerate}
\end{rem}

 The next lemma provides an estimate for the $H^k$-norm of
the solution $({\varphi},{\theta},{\psi},{\pi})$ with $\xi$ varying,
which will be useful in the next section when such a solution is
integrated in a Fourier synthesis. To emphasize the dependence on
$\xi$, we write these solutions as
$({\varphi}(\xi)={\varphi}(\xi,x_3),{\theta}(\xi)
={\theta}(\xi,x_3),{\psi}(\xi)={\psi}(|\xi|,x_3),{\pi}(\xi)
={\pi}(|\xi|,x_3))$.
\begin{lem}\label{lem:0201}
Let $\xi\in \mathbb{R}^2$ with $0< R_1<|\xi|<R_2$, $\varphi(\xi)$,
$\theta(\xi)$, $\psi(\xi)$, $\pi(\xi)$ and $\lambda(\xi)$ be
constructed as in Theorem \ref{thm:0302}, then for any $k\geq  0$
there exit positive constants $A_k$, $B_k$, $C_k$  and  $D$, which may
depend on $R_1$, $R_2$, $\bar{\rho}$, $\mu$ and $g$, such that
\begin{eqnarray}
&&\label{0245} \|\psi(\xi)\|_{H^k(\mathbb{R})}\leq A_k,\\
&&\label{0246} \|\pi(\xi)\|_{H^k(\mathbb{R})}\leq B_k,\\
&&\label{0247}\|\varphi(\xi)\|_{H^k(\mathbb{R})}+\|\theta(\xi)\|_{H^k(\mathbb{R})}
\leq C_k.
\end{eqnarray} Moreover,
\begin{equation}\label{0248}
\|\psi\|_{L^2(\mathbb{R})}^2> 0.
\end{equation}
\end{lem}
\begin{pf} Throughout this proof, we denote by $\tilde{c}$ a generic positive constant which may
vary from line to line, and depend on $R_1$, $R_2$, $\bar{\rho}$, $\mu$ and $g$.

(i) First, since  $\psi\in \mathcal{A}$, we see that  (\ref{0248})
holds, and there exists a constant $\tilde{c}$, such that
\begin{equation}\label{0249}
\|\psi(\xi)\|_{H^1(\mathbb{R})}\leq \tilde{c}.
\end{equation}  On the other hand, in view of Proposition
\ref{pro:0205}, we have
\begin{equation}\label{0250}
 \lambda(\xi) \geq \tilde{c}>0\quad\mbox{for any }|\xi|\in (R_1,R_2).
\end{equation}
Similarly to (\ref{0236}), we use (\ref{0249}), (\ref{0250}) and (\ref{0235}) with $\lambda(\xi)$
in place of $\alpha(\xi)$ to deduce that
\begin{equation}\label{0251}
\|\psi(\xi)\|_{H^2(\mathbb{R})}\leq \tilde{c}.
\end{equation}

 We now rewrite (\ref{0207}) as
\begin{equation}\label{0252}
\psi''''(\xi)=\big[{\lambda(\lambda\bar{\rho}+2\mu
|\xi|^2)\psi''(\xi)+\lambda^2\bar{\rho}'\psi'(\xi)-|\xi|^2(\lambda^2\bar{\rho}
+\lambda\mu |\xi|^2-g\bar{\rho}')\psi(\xi)}\big]/\lambda\mu,
\end{equation}
which, together with (\ref{0221}), (\ref{0250}) and (\ref{0251}), yields
\begin{equation}\label{0253}
\|\psi(\xi)\|_{H^4(\mathbb{R})}\leq \tilde{c}.
\end{equation}

 Differentiating (\ref{0252}) with respect to $x_3$ and
using (\ref{0253}), we find, by induction on $k$, that (\ref{0245})
holds for any $k\geq 0$. Consequently, we can employ (\ref{0245})
and the expression (\ref{0243}) of $\pi$ to deduce that (\ref{0246})
holds for any $k\geq 0$.

(ii) Finally, we verify (\ref{0247}). By virtue of (\ref{0244}),
$(\varphi,\theta)=-\psi'(\xi_1,\xi_2)/|\xi|^2$. Hence, from (\ref{0249}) we get
\begin{equation}\label{0254}\|\varphi(\xi)\|_{L^2(\mathbb{R})}
+\|\theta(\xi)\|_{L^2(\mathbb{R})}<\tilde{c}.\end{equation}
Noticing that (\ref{0205})$_1$ and (\ref{0205})$_2$ can be rewritten as
\begin{equation}\label{0255}
\varphi''(\xi)=\left(\frac{\lambda\bar{\rho} }{\mu}+ |\xi|^2\right)\varphi(\xi)-\frac{\xi_1 \pi(\xi)}{\mu},\\
\end{equation}
and
\begin{equation}
\label{0256}\theta''(\xi)=\left(\frac{\lambda\bar{\rho} }{\mu}+
|\xi|^2\right)\theta(\xi)-\frac{\xi_2 \pi(\xi)}{\mu}, \end{equation}
we apply (\ref{0254}) to (\ref{0255}) and (\ref{0256}) to obtain
\begin{equation}\label{0257}\begin{aligned}
\|\varphi{''}(\xi)\|_{L^2(\mathbb{R})}+\|\theta{''}(\xi)\|_{L^2(\mathbb{R})}\leq&
\tilde{c}.
\end{aligned}\end{equation}

Using (\ref{0254}) and (\ref{0257}), analogously to (\ref{0221}),
we infer that
\begin{equation}\label{0258}\begin{aligned}
\|\varphi{'}(\xi)\|_{L^2(\mathbb{R})}+\|\theta{'}(\xi)\|_{L^2(\mathbb{R})}\leq&
\tilde{c}.
\end{aligned}\end{equation}
Putting  (\ref{0254})--(\ref{0258}) together, we immediately
obtain (\ref{0247}). This completes the proof. \hfill $\Box$
\end{pf}

\subsection{Exponential growth rate}

In this subsection we use the Fourier synthesis to build growing
solutions to (\ref{0108}) out of the solutions constructed in the
previous subsection (Theorem \ref{thm:0302}) for any fixed spatial
frequency $\xi\in \mathbb{R}^2$ with $|\xi|>0$. The constructed
solutions will grow in-time in the Sobolev space of order $k$.
\begin{thm}\label{thm:0203}
 Let $0<R_1<R_2<\infty$ and $f\in C_0^\infty(R_1,R_2)$ be a real-valued
function. For $\xi\in \mathbb{R}^2$ with $ |\xi|\in (0, \infty)$,
define
\begin{equation*} \mathbf{v}(\xi,
x_3)=-i\varphi(\xi,x_3)e_1-i\theta(\xi,x_3)e_2+\psi(\xi, x_3)e_3,
\end{equation*}
where $(\varphi, \theta, \psi, \pi)(\xi, x_3)$ with
$\lambda(|\xi|)>0$ is the solution given by Theorem \ref{thm:0302}.
Denote
\begin{eqnarray}
&&\label{0259}\varrho(t,\mathbf{x})=-\frac{\bar{\rho}'(x_3)}{4\pi^2}\int_{\mathbb{R}^2}
f(|\xi|){v}_3(\xi,x_3)e^{\lambda(|\xi|)t}e^{ix'\xi}\mathrm{d}\xi ,\\
&&\label{0260}\mathbf{u}(t,\mathbf{x})=\frac{1}{4\pi^2}\int_{\mathbb{R}^2}
\lambda(\xi)f(|\xi|)\mathbf{v}(\xi,x_3)e^{\lambda(|\xi|)t}e^{ix'\xi}\mathrm{d}\xi ,\\
&&\label{0261}q(t,\mathbf{x})=\frac{1}{4\pi^2}\int_{\mathbb{R}^2}
\lambda(\xi)f(|\xi|){\pi}(\xi,x_3)e^{\lambda(|\xi|)t}e^{ix'\xi}\mathrm{d}\xi,
\end{eqnarray}
Then, $(\varrho,\mathbf{u},q)$ is a real-valued solution to the linearized problem
 (\ref{0108}) along with (\ref{0107}). For every $k\in \mathbb{N}$, we have the estimate
\begin{equation}\label{0262}
\|\varrho(0)\|_{H^k}+\|\mathbf{u}(0)\|_{H^k} +\|q(0)\|_{H^k} \leq {D}_k
\left(\int_{\mathbb{R}^2}(1+|\xi|^2)^{k+2}|f(|\xi|)|^2\mathrm{d}\xi\right)^{1/2}<\infty ,
\end{equation}
where $D_k>0$ is a constant depending on $k$,
$\bar{\rho}$, $R_1$, $R_2$ and $g$. Moreover, for every $t>0$ we
have $(\varrho(t),\mathbf{u}(t),q(t))\in H^k$, and
\begin{eqnarray}
&&\label{0263} e^{t\lambda_0(f)}\|\varrho(0)\|_{H^k}
\leq \|\varrho(t)\|_{H^k}\leq e^{t\Lambda}\|\varrho(0)\|_{H^k},\\[1mm]
&&\label{0264}e^{t\lambda_0(f)}\|u_i(0)\|_{H^k}\leq \|u_i(t)\|_{H^k}\leq
e^{t\Lambda}\|u_i(0)\|_{H^k},\quad i=1,2,3,\\[1mm]
&&\label{0265}e^{t\lambda_0(f)}\|q(0)\|_{H^k}\leq
\|q(t)\|_{H^k}\leq e^{t\Lambda}\|q(0)\|_{H^k},
\end{eqnarray} where
\begin{equation}\label{0266}\lambda_0(f)=\inf_{|\xi|\in
\mathrm{supp}(f)}\lambda(|\xi|)>0\end{equation} and
 \begin{equation}\label{0267}\Lambda=\sup_{0<|\xi|<+\infty}\lambda(|\xi|)<\sqrt{g}
 \left\|\sqrt{{\bar{\rho}'}/{\bar{\rho}}}
 \right\|_{L^\infty(\mathbb{R})}.\end{equation}
 In particular, \begin{equation}\label{0268}
\|{u}_3(0)\|_{H^k}>0\;\;\mbox{ if }f\not\equiv 0,\end{equation}
and we can further take proper constants $R_1$, $R_2$, such that
\begin{equation}\label{0269}\lambda_0(f)=\Lambda/2.\end{equation}
\end{thm}
\begin{pf} Obviously, (\ref{0266}), (\ref{0267}) and (\ref{0269}) follow from Proposition \ref{pro:0205}.
For each fixed $\xi\in \mathbb{R}^2$,
\begin{eqnarray*}
&& \tilde{\varrho} (t,\textbf{x})=-\bar{\rho}'f(|\xi|){v_3}(\xi,x_3)e^{\lambda(|\xi|)t}e^{ix'\xi},\\[1mm]
&& \tilde{\mathbf{u}}(t,\textbf{x})=\lambda(|\xi|)f(|\xi|)\mathbf{v}(\xi,x_3)e^{\lambda(|\xi|)t}e^{ix'\xi},\\[1mm]
&&
\tilde{q}(t,\textbf{x})=\lambda(|\xi|)f(|\xi|){\pi}(\xi,x_3)e^{\lambda(|\xi|)t}e^{ix'\xi}
\end{eqnarray*}
give a solution to (\ref{0108}). Since $f\in C_0^\infty(R_1,R_2)$,
Lemma \ref{lem:0201} implies that
\begin{equation*}
\sup_{\xi\in\mathrm{supp}(f)}\|\partial_{3}^k(\tilde{\varrho},\tilde{\mathbf{u}},\tilde{q})(\xi,\cdot)\|_{L^\infty}<\infty\quad
\mbox{ for all }k\in \mathbb{N}.\end{equation*}
 These bounds show that the Fourier
synthesis of the solution given by (\ref{0259})--(\ref{0261}) is
also a solution of (\ref{0108}). Because $f$ is real-valued and
radial, by Remark \ref{rem:0201} we can easily verify
that the Fourier synthesis is real-valued.

The estimates (\ref{0262}) and (\ref{0268}) follow from Lemma
\ref{lem:0201} with arbitrary $k\geq 0$, and the fact that $f$ is
compactly supported and
\begin{equation}\label{0270}\begin{aligned}
\|r\|^2_{{{H}}^k}:=&\sum_{j=0}^k\int_{\mathbb{R}^2\times
\mathbb{R}}(1+|\xi|^2)^{k-j}\left|\partial_{x_3}^j\hat{r}(\xi,x_3)\right|^2\mathrm{d}\xi\mathrm{d}x_3\\
=
&\sum_{j=0}^k\int_{\mathbb{R}^2}(1+|\xi|^2)^{k-j}\left\|\partial_{x_3}^j \hat{r}(\xi,\cdot)
\right\|^2_{L^2(\mathbb{R})}\mathrm{d}\xi\quad\mbox{for }r=\varrho,\mbox{ or }u_i,\mbox{ or }q.
\end{aligned}\end{equation}
Finally, we can use (\ref{0266}),
(\ref{0267}) and (\ref{0270}) to obtain the estimates
(\ref{0263})--(\ref{0265}).
 \hfill $\Box$
\end{pf}

\section{Uniqueness of the linearized equations}
In this section, we will show the uniqueness of solutions to the
linearized problem, which will be used in the proof of Theorem
\ref{thm:0101} in Section 5. We first define the function space of strong solutions.
\begin{equation*}\begin{aligned}
\mathcal{Q}(T):=\{(\varrho,\mathbf{u},q)~|&~\varrho\in
C^0([0,T],L^2),\ \nabla q\in L^2(0,T; H^1_{\mathrm{loc}}),\\
&~ \mathbf{u}\in C^0([0,T],(L^2)^3)\cap L^2(0,T;(H^2)^3),\;
\partial_t\mathbf{u}\in (L^2((0,T)\times\mathbb{R}^3))^3\}.\end{aligned}
\end{equation*}
We claim that the solution to the linearized problem is unique in the function space $\mathcal{Q}$.
\begin{thm} {\rm (Uniqueness)} \label{thm:0301}
Assume that $(\tilde{\varrho}, \tilde{\mathbf{u}}, \tilde{q})$,
$(\bar{\varrho}, \bar{\mathbf{u}}, \bar{q})\in\mathcal{Q}(T)$ are two strong solutions of
(\ref{0108}) with
$(\tilde{\varrho},\tilde{\mathbf{u}})(0)=(\bar{\varrho},\bar{\mathbf{u}})(0)$.
Then, $(\tilde{\varrho},\tilde{\mathbf{u}},\nabla
\tilde{q})=(\bar{\varrho},\bar{\mathbf{u}},\nabla \bar{q})$.
\end{thm}
\begin{pf}
Let $(\varrho,\mathbf{u},
q)=(\tilde{\varrho}-\bar{\varrho},\tilde{\mathbf{u}}-\bar{\mathbf{u}},
\tilde{q}-\bar{q})$. Then $(\varrho,\mathbf{u}, q)\in
\mathcal{Q}(T)$ is still a strong solution to the linearized system
(\ref{0108}) with $(\varrho ,\mathbf{u})(0)=(0,\mathbf{0})$.

Multiplying (\ref{0108})$_2$ by
$\mathbf{\varphi}\in (C_0^\infty((0,t)\times \mathbb{R}^3))^3$ with
$\mathrm{div}\varphi=0$, integrating over $(0,t)\times\mathbb{R}^3$, and
then integrating by parts, we arrive at
\begin{equation*}\label{0301}\begin{aligned}
&\int_0^{t}\int_{\mathbb{R}^3} \bar{\rho}\partial_t\mathbf{u}\cdot
\varphi\mathrm{d}\mathbf{x}\mathrm{d}\tau
+\int_0^t\int_{\mathbb{R}^3}\mu\nabla \mathbf{u}:\nabla
\varphi\mathrm{d}\mathbf{x}\mathrm{d}\tau
=-g\int_0^{t}\int_{\mathbb{R}^3} \rho
\varphi_3\mathrm{d}\mathbf{x}\mathrm{d}\tau.
\end{aligned}\end{equation*}
Using the standard density argument (see, for instance, \cite{RDJLLA} or \cite[Section 2.1]{LPLMTFM96}),
we could take $\mathbf{u}$ as a test function to get
\begin{equation}\label{0302}\begin{aligned}
&\int_0^{t}\int_{\mathbb{R}^3} \bar{\rho}\partial_t\mathbf{u}\cdot
\mathbf{u}\mathrm{d}\mathbf{x}\mathrm{d}\tau
+\int_0^t\int_{\mathbb{R}^3}\mu\nabla \mathbf{u}:\nabla
\mathbf{u}\mathrm{d}\mathbf{x}\mathrm{d}\tau
=-g\int_0^{t}\int_{\mathbb{R}^3} \varrho
u_3\mathrm{d}\mathbf{x}\mathrm{d}\tau ,
\end{aligned}\end{equation}
which will be employed to derive an energy-like estimate. In fact, recalling
$$\mathbf{u}\in L^2(0,T;(H^1)^3)\cap C^0([0,T],(L^2)^3)\mbox{ and }
\partial_t\mathbf{u}\in (L^2((0,T)\times \mathbb{R}^3))^3,$$
and $\mathbf{u}(0)=\mathbf{0}$, we easily deduce
\begin{equation}\label{0303}\begin{aligned}
\int_0^{t}\int_\Omega \bar{\rho}
\partial_t\mathbf{u}\cdot \mathbf{u}\mathrm{d}\mathbf{x}\mathrm{d}\tau =\frac{1}{2}\int_\Omega
\bar{\rho} \mathbf{u}^2(t)\mathrm{d}\mathbf{x}.
\end{aligned}\end{equation}
 Since $\varrho\in C^0([0,T],L^2)$ and $\varrho(0)=0$, the equation (\ref{0108})$_1$ gives
\begin{equation*}\label{0304}
\varrho(t,\mathbf{x})=\int_0^t\bar{\rho}'{u}_3(s,\mathbf{x})\mathrm{d}s\mbox{
for any }t\geq 0,
\end{equation*}
 Consequently, with the help of the regularity of $\partial_iu_3$, the
property of absolutely continuous functions and Fubini's theorem, we conclude that
\begin{equation}\label{0305}\begin{aligned}
\int_0^{t}\int_{\mathbb{R}^3} \varrho
u_3\mathrm{d}\mathbf{x}\mathrm{d}\tau
=&\int_{\mathbb{R}^3}\bar{\rho}'\int_0^{t}\int_0^\tau{u}_3(s,\mathbf{x})\mathrm{d}s
u_3(\tau,\mathbf{x})\mathrm{d}\tau\mathrm{d}\mathbf{x}\\
=&\frac{1}{2}\int_{\mathbb{R}^3}\bar{\rho}'\int_0^{t}\frac{d}{d\tau}\left(\int_0^\tau{u}_3(s,\mathbf{x})\mathrm{d}s
\right)^2\mathrm{d}\tau\mathrm{d}\mathbf{x}\\
=&\frac{1}{2}\int_{\mathbb{R}^3}\bar{\rho}'\left(\int_0^{t}{u}_3(\tau,\mathbf{x})\mathrm{d}\tau
\right)^2\mathrm{d}\mathbf{x}\\
 \leq &\frac{t}{2}\int_0^{t}\int_{\mathbb{R}^3}\bar{\rho}'{u}_3^2\mathrm{d}\mathbf{x}\mathrm{d}\tau\\
 \leq & \frac{t}{2}\left\|\frac{\bar{\rho}'}{\bar{\rho}}\right\|_{L^\infty(\mathbb{R})}
 \int_0^{t}\int_{\mathbb{R}^3}\bar{\rho}\mathbf{u}^2
\mathrm{d}\mathbf{x}\mathrm{d}\tau .
\end{aligned}\end{equation}

Substituting (\ref{0303}) and (\ref{0305}) into (\ref{0302}),  we conclude that
\begin{equation*}\label{0306}\begin{aligned}
&\frac{1}{2}\int_{\mathbb{R}^3} \bar{\rho}
\mathbf{u}^2(t)\mathrm{d}\mathbf{x}
+\int_0^t\int_{\mathbb{R}^3}\mu\nabla \mathbf{u}:\nabla
\mathbf{u}\mathrm{d}\mathbf{x}\mathrm{d}\tau \leq
\frac{t}{2}\left\|\frac{\bar{\rho}'}{\bar{\rho}}\right\|_{L^\infty(\mathbb{R}^3
 )}\int_0^{t}\int_{\mathbb{R}^3}\bar{\rho}\mathbf{u}^2
\mathrm{d}\mathbf{x}\mathrm{d}\tau,
\end{aligned}\end{equation*}
which yields
\begin{equation}\label{0307}\|\sqrt{ \bar{\rho}}\mathbf{u}(t)\|^2_{L^2}\leq
T\left\|\frac{\bar{\rho}'}{\bar{\rho}}\right\|_{L^\infty(\mathbb{R})}
\int_0^{t}\|\sqrt{\bar{\rho}}\mathbf{u}\|_{L^2}\mathrm{d}\tau.
\end{equation}

Applying Grownwall's inequality to (\ref{0307}), we get
\begin{equation*}
\|\sqrt{
\bar{\rho}}\mathbf{u}(t)\|^2_{L^2}=0\quad\mbox{ for any }t\in [0,T],
\end{equation*}
which yields $\mathbf{u}=0$, i.e.,
$\tilde{\mathbf{u}}=\bar{\mathbf{u}}$, since $\bar{\rho}>0$. This,
combined with (\ref{0108})$_1$ and (\ref{0108})$_2$, proves that
\begin{equation*}  (\tilde{\rho},\tilde{\mathbf{u}},\nabla
\tilde{q})=(\bar{\rho},\bar{\mathbf{u}},\nabla \bar{q})\quad\mbox{ for any }t\in (0,T].\end{equation*}
Thus, the desired conclusion follows.   \hfill $\Box$
\end{pf}

\section{Nonlinear energy estimates of the perturbed problem}

In this section, we shall derive some  nonlinear energy estimates
for the perturbed problem, which will also be used in the proof of
Theorem \ref{thm:0101} in Section 5. To this end, let
$(\varrho,\mathbf{u},q)$ be a classical solution of the perturbed
problem (\ref{0105})--(\ref{0107}) with $\rho:=\varrho+\bar{\rho}>0$
in $[0,T]\times \mathbb{R}^3$ for some $T>0$. Moreover, we assume
that the classical solution $(\varrho,\mathbf{u},q)$ satisfies the
initial condition
\begin{equation}\label{0401}
\sqrt{\|\varrho_0\|_{H^1}^2+\|\mathbf{u}_0\|_{H^2}^2}=\delta_0\leq
1.\end{equation} The restricted relation between $T$ and $\delta_0$
will be given at the end of Subsection 4.1.

In what follows, we denote by $C_i$ ($i=1,2,\cdots$), $C(T)$ and $C$
generic positive constants depending on $\mu$, $g$ and $\bar{\rho}$.
In addition, $C(T)$ also depends on $T$ and is nondecreasing with
respect to $T$; the subscript $i$ of $C_i$ emphases that we may
repeat to  use the constant $C_i$ in the process of estimates.

\subsection{Estimates  for $\|\mathbf{u}_t\|_{L^2}$ and $\|\nabla
\mathbf{u}\|_{H^1}$}
 We first observe that the continuity equation
(\ref{0105})$_1$ and the incompressibility condition
(\ref{0105})$_3$ imply immediately that
 for any
$t\in (0,T]$,
\begin{equation}\label{0402}\alpha:=\inf_{\mathbf{x}\in\mathbb{R}^3}\{\rho_0(\mathbf{x})\}
\leq \rho(t)\leq \sup_{\mathbf{x}\in
\mathbb{R}^3}\{\rho_0({\mathbf{x}})\}:=\beta\mbox{ or
}\alpha-\bar{\rho}\leq \varrho(t)\leq \beta-\bar{\rho},
\end{equation}
and
\begin{equation}\label{0403}\frac{d}{dt}
\|\varrho(t)\|_{L^2}^2=-2\int_{\mathbb{R}^3} \bar{\rho}'\varrho
u_3\mathrm{d}\mathbf{x}\leq2
\left\|\frac{\bar{\rho}'\varrho}{\sqrt{{\rho}}}\right\|_{L^2}\|\sqrt{{\rho}}\mathbf{u}\|_{L^2}\leq
2
\alpha^{-1}\|\bar{\rho}'\|_{L^\infty}\left\|\varrho\right\|_{L^2}\|\sqrt{{\rho}}\mathbf{u}\|_{L^2}.
\end{equation}

Multiplying (\ref{0105})$_2$ by $\mf{u}$, using (\ref{0105})$_1$, and then
integrating (by parts) over $(0,t)\times \mathbb{R}^3$, we obtain
\begin{equation*}\label{0404}\frac{1}{2}\frac{d}{dt}\int_{\mathbb{R}^3}
\rho|\mathbf{u}|^2(t)\mathrm{d}\mathbf{x}+\mu\int_{\mathbb{R}^3}|\nabla \mathbf{u}|^2
\mathrm{d}\mathbf{x}= -g\int_{\mathbb{R}^3}{ \varrho
}{u}_3\mathrm{d}\mathbf{x}.\end{equation*}
Since the integral on the right-hand side is bounded from above by
$g\alpha^{-\frac{1}{2}}\|\varrho\|_{L^2}\|\sqrt{\rho}\mathbf{u}\|_{L^2}$, we get
\begin{equation}\label{0405}\frac{d}{dt}\|\sqrt{\rho}\mathbf{u}(t)\|_{L^2}^2+\mu\|\nabla
\mathbf{u}\|_{L^2}^2\leq
2g\alpha^{-\frac{1}{2}}\left\|\varrho\right\|_{L^2}\|\sqrt{{\rho}}\mathbf{u}\|_{L^2}.\end{equation}
Combining (\ref{0403}) with (\ref{0405}) and using Cauchy-Schwarz's inequality, we obtain
\begin{equation}\label{0406}\begin{aligned}
\frac{d}{dt}(\|\varrho(t)\|_{L^2}^2+\|\sqrt{{\rho}}\mathbf{u}(t)\|_{L^2}^2)+\mu\|\nabla
\mathbf{u}\|_{L^2}^2 \leq &
C_1(\|\varrho(t)\|_{L^2}^2+\|\sqrt{{\rho}}\mathbf{u}(t)\|_{L^2}^2),\end{aligned}\end{equation}
which implies
\begin{equation}\label{0407}\|\varrho(t)\|_{L^2}^2+\|\sqrt{{\rho}}\mathbf{u}(t)\|_{L^2}^2\leq
\delta_0e^{C_1t}.
\end{equation}
 In particular, making use of (\ref{0402}), (\ref{0406}) and (\ref{0407}), we arrive at
\begin{equation}\label{0408}\|\varrho(t)\|_{L^2}^2+\|\mathbf{u}(t)\|_{L^2}^2
+\int_0^t\|\nabla \mathbf{u}(s)\|_{L^2}^2\mm{d}\tau\leq C\delta_0e^{C_1t}.
\end{equation}

To control $\mathbf{u}_t$, we multiply (\ref{0105})$_2$ by $\mathbf{u}_t$ in $L^2$ and
apply Cauchy-Schwarz's inequality to infer that
\begin{equation}\label{0410} \frac{1}{2}\|\sqrt{\rho}\mathbf{u}_t\|_{L^2}^2+\mu\frac{d}{dt}
\|\nabla \mathbf{u}(t)\|_{L^2}^2\leq C(\|{ \varrho}\|_{L^2}^2
+\|\sqrt{\rho}\mathbf{u}\cdot\nabla
\mathbf{u}\|_{L^2}^2).\end{equation}
To bound the second term on the right-hand side of (\ref{0410}), we recall
that $(\mathbf{u},q)$ is a solution of the Stokes equations:
\begin{equation*}\label{0411}
-\mu \Delta\mathbf{u} +\nabla q = -\rho \mathbf{u}_t
-\rho (\mathbf{u}\cdot \nabla \mathbf{u})-\varrho {g}e_3,\quad
\mathrm{div}\,\mathbf{u}=0 \mbox{ in }\mathbb{R}^3.
\end{equation*}
Therefore, the classical regularity theory on the Stokes equations (see \cite[Theorem 2.1]{GGPA94}) gives
\begin{equation*} 
\begin{aligned}\|\nabla^2 \mathbf{u} \|_{L^2}^2\leq &
C_2(\|\sqrt{\rho}\mathbf{u}_t\|^2_{L^2}+\|\sqrt{\rho}\mathbf{u}\cdot\nabla
\mathbf{u}\|_{L^2}^2+g\|{\varrho}\|_{L^2}^2),
\end{aligned}\end{equation*}
which combined with (\ref{0410}) results in
\begin{equation}\label{0413}
\frac{1}{4}
\|\sqrt{\rho}\mathbf{u}_t\|_{L^2}^2+\varepsilon\|\nabla^2
\mathbf{u}\|_{L^2}^2+\mu\frac{d}{dt} \|\nabla
\mathbf{u}(t)\|_{L^2}^2\leq C_3(\|{ \varrho}\|_{L^2}^2
+\|\mathbf{u}\cdot\nabla\mathbf{u}\|_{L^2}^2)\;\;\mbox{ with }C_3>1.
\end{equation}

If we apply H\"{o}lder's and Sobolev's inequalities, we have
\begin{equation*}\label{0414}\begin{aligned}
\|\mathbf{u}\cdot\nabla\mathbf{u}\|_{L^2}^2\leq &
\|\mathbf{u}\|_{L^6}^2\|\nabla \mathbf{u}\|_{L^3}^2\leq
\|\mathbf{u}\|_{L^6}^2\|\nabla \mathbf{u}\|_{L^2}\|\nabla
\mathbf{u}\|_{L^6}\\
\leq & C_0\|\nabla \mathbf{u}\|_{L^2}^3\|\nabla
\mathbf{u}\|_{H^1}\leq \frac{2C^2_0 C_3}{\varepsilon}\|\nabla
\mathbf{u}\|_{L^2}^6+\frac{\varepsilon}{8 C_3}\|\nabla
\mathbf{u}\|_{H^1}^2.
\end{aligned}\end{equation*}
Substituting the above inequality into (\ref{0413}) and integrating over
($0,t$), we conclude that
\begin{equation}\label{n0411}
\frac{1}{4}
\|\sqrt{\rho}\mathbf{u}_t\|_{L^2}^2+\frac{7\varepsilon}{8}\|\nabla^2
\mathbf{u}\|_{L^2}^2+\mu\frac{d}{dt} \|\nabla
\mathbf{u}(t)\|_{L^2}^2\leq C(\| \varrho\|_{L^2}^2 +\|\nabla
\mathbf{u}\|^6_{L^2}).
\end{equation}

 Letting $\delta_0 e^{C_1T}\leq 1$, we get from (\ref{0406}) and (\ref{0413}) that
\begin{equation}\label{0416}  \begin{aligned}
\mu \|\nabla \mathbf{u}(t)\|_{L^2}^2\leq & \int_0^t
C(\|\varrho\|_{L^2}^2 +\|\nabla \mathbf{u}\|^6_{L^2})\mathrm{d}s+\mu
\|\nabla \mathbf{u}_0\|_{L^2}^2  \\
\leq & C\delta_0^2(e^{C_1t}-1) +C\int_0^t\|\nabla
\mathbf{u}\|^6_{L^2}\mathrm{d}s+\mu \|\nabla \mathbf{u}_0\|_{L^2}^2\\
\leq & C_4\Big(\delta_0+\int_0^t\|\nabla \mathbf{u}\|^6_{L^2}\mathrm{d}s\Big),
\end{aligned} \end{equation}
which yields
\begin{equation}\label{0417}
\|\nabla \mathbf{u}(t)\|_{L^2}^2\leq
\sqrt{\frac{\delta_0^2C_4^2}{1-2t\delta_0^2C_4^3}}.
\end{equation}

Now, we take
\begin{equation*}
T=\min\left\{\frac{1}{C_1}\mathrm{ln}\frac{1}{\delta_0},\frac{1}{4\delta_0^2C_4^3}\right\}.
\end{equation*}
In particular, there exists a sufficiently small constant
$\delta_1>0$, such that
\begin{equation}\label{0419}
T= \frac{-\mathrm{ln}\delta_0}{C_1}\quad\mbox{ for any }\delta_0\in
(0,\delta_1).
\end{equation}

\emph{From now on, we always assume that $\delta_0$ and $T$ satisfy the relation (\ref{0419}).}
Thus, (\ref{0417}) gives
\begin{equation}\label{0420}
\|\nabla \mathbf{u}(t)\|_{L^2}^2\leq \sqrt{2{\delta_0^2C_4^2}}\leq
C\delta_0\quad\mbox{ for any }t\in (0,T].
\end{equation}
Making use of (\ref{0402}), (\ref{0407}), ({\ref{0420}) and the
first inequality in (\ref{0416}), we deduce from (\ref{n0411}) that
\begin{equation}\label{0421}\|\varrho(t)\|_{L^2}^2+\|\mathbf{u}(t)\|_{H^1}^2+\int_0^t(\|
\mathbf{u}_t(s)\|_{L^2}^2+\|\nabla
\mathbf{u}(s)\|_{H^1}^2)\mm{d}s\leq C(T)\delta_0^2\mbox{ for any
}t\in (0,T].
\end{equation}

\subsection{Estimates for $\|\nabla \mathbf{u}_t\|_{L^2}$ and $\|\nabla^2\mathbf{u}\|_{L^6}$}
Using (\ref{0105})$_1$ and keeping in mind that $p=q+\bar{p}$, we can rewrite (\ref{0105})$_2$ as
\begin{equation*}\label{}
\rho\mathbf{u}_t+\rho\mathbf{u}\cdot \nabla \mathbf{u}+\nabla
p=\mu\Delta\mathbf{u}-\rho g e_3,
\end{equation*}
whence, by taking the time derivative,
\begin{equation*}\label{}
\rho\mathbf{u}_{tt}+\rho\mathbf{u}\cdot \nabla
\mathbf{u}_t-\mu\Delta\mathbf{u}_t+\nabla
p_t=-\rho_t(\mathbf{u}_t+\mathbf{u}\cdot\nabla \mathbf{u}+g
e_3)-\rho\mathbf{u}_t\cdot\nabla \mathbf{u},
\end{equation*}
which, by using the continuity equation, can be written as
\begin{equation*}\label{}\begin{aligned}
&\rho\left(\frac{1}{2}|\mathbf{u}_{t}|^2\right)_t+\rho\mathbf{u}\cdot
\nabla
\left(\frac{1}{2}|\mathbf{u}_t|^2\right)-\mu\Delta\mathbf{u}_t\cdot\mathbf{u}_t+\nabla
p_t\cdot\mathbf{u}_t\\
&=\mathrm{div}(\rho\mathbf{u})(\mathbf{u}_t+\mathbf{u}\cdot\nabla
\mathbf{u}+ g e_3)\cdot\mathbf{u}_t-\rho(\mathbf{u}_t\cdot\nabla
\mathbf{u})\mathbf{u}_t.\end{aligned}
\end{equation*}
Hence, by integrating by parts, we see that
\begin{equation}\label{0425}\begin{aligned}
&\frac{1}{2}\frac{d}{dt}\int_{\mathbb{R}^3}\rho|\mathbf{u}_{t}(t)|^2\mathrm{d}\mathbf{x}
+\mu\int_{\mathbb{R}^3}|\nabla \mathbf{u}_t|^2\mathrm{d}\mathbf{x}\\
&\leq \int_{\mathbb{R}^3}2\rho|\mathbf{u}||\mathbf{u}_t||\nabla
\mathbf{u}_t|+\rho|\mathbf{u}||\mathbf{u}_t||\nabla
\mathbf{u}|^2+\rho |\mathbf{u}|^2|\mathbf{u}_t||\nabla^2 \mathbf{u}|\\
&\quad +\rho|\mathbf{u} |^2|\nabla \mathbf{u}||\nabla
\mathbf{u}_t|+\rho|\mathbf{u}_t|^2|\nabla
\mathbf{u}|+g\rho|\mathbf{u}||\nabla
\mathbf{u}_t|:=\sum_{i=1}^6I_j,\end{aligned}
\end{equation}
where $I_j$ can be bounded as follows, employing straightforward calculations.
\begin{eqnarray*}
&& I_1\leq
2\|\rho\|_{L^\infty}^{1/2}\|\mathbf{u}\|_{L^6}\|\sqrt{\rho}\mathbf{u}_t\|_{L^3}\|\nabla
\mathbf{u}_t\|_{L^2}\leq 2
\|\rho\|_{L^\infty}^{1/2}\|\mathbf{u}\|_{L^6}\|\sqrt{\rho}\mathbf{u}_t\|_{L^2}^{1/2}
\|\sqrt{\rho}\mathbf{u}_t\|_{L^6}^{1/2}\|\nabla\mathbf{u}_t\|_{L^2}\\
&& \quad \leq C\|\rho\|_{L^\infty}^{3/4}\|\nabla
\mathbf{u}\|_{L^2}\|\sqrt{\rho}\mathbf{u}_t\|_{L^2}^{1/2}\|\nabla
\mathbf{u}_t\|_{L^2}^{3/2}\leq C(\varepsilon)\|\nabla
\mathbf{u}\|_{L^2}^4\|\sqrt{\rho}\mathbf{u}_t\|_{L^2}^2+\varepsilon\|\nabla
\mathbf{u}_t\|_{L^2}^2,\\[1mm]
&& I_2\leq C\|\rho\|_{L^\infty}\|\nabla \mathbf{u} \|_{L^2}\|\nabla
\mathbf{u}_t\|_{L^2}\|\nabla \mathbf{u}\|_{L^2}\|\nabla
\mathbf{u}\|_{H^1}\leq C(\varepsilon)\|\nabla
\mathbf{u}\|_{L^2}^4\|\nabla
\mathbf{u}\|_{H^1}^2+\varepsilon\|\nabla\mathbf{u}_t \|_{L^2}^2,\\[1mm]
&&I_3\leq C\|\rho\|_{L^\infty}\|\nabla \mathbf{u}\|_{L^2}^2\|\nabla
\mathbf{u}_t\|_{L^2}\|\nabla^2 \mathbf{u}\|_{L^2}\leq
C(\varepsilon)\|\nabla \mathbf{u}\|_{L^2}^4\|\nabla^2
\mathbf{u}\|_{L^2}^2+\varepsilon\|\nabla \mathbf{u}_t\|_{L^2}^2,\\[1mm]
&& I_4\leq C\|\rho\|_{L^\infty}\|\nabla \mathbf{u}\|_{L^2}^2\|\nabla
\mathbf{u}\|_{H^1}\|\nabla \mathbf{u}_t\|_{L^2}\leq
C(\varepsilon)\|\nabla \mathbf{u}\|_{L^2}^4\|\nabla
\mathbf{u}\|_{H^1}^2+\varepsilon\|\nabla \mathbf{u}_t\|_{L^2},\\
&& I_5\leq C\|\rho\|_{L^\infty}^{3/4}\|\nabla
\mathbf{u}\|_{L^2}\|\sqrt{\rho}\mathbf{u}_t\|_{L^2}^{1/2}\|\nabla
\mathbf{u}_t\|_{L^2}^{3/2}\leq C(\varepsilon)\|\nabla
\mathbf{u}\|_{L^2}^4\|\sqrt{\rho}\mathbf{u}_t\|_{L^2}^2+\varepsilon\|\nabla
\mathbf{u}_t\|_{L^2}^2, \\[1mm]
&& I_6\leq C(\varepsilon)\|\mathbf{u}\|_{L^2}^2+\varepsilon\|\nabla\mathbf{u}_t\|_{L^2}^2,
\end{eqnarray*}
where $C(\varepsilon)$ is a positive constant which may depend on $\varepsilon$.
Inserting all the above estimates into (\ref{0425}), we conclude
\begin{equation*}\label{0427}
\frac{d}{dt}\|\sqrt{\rho}\mathbf{u}_{t}(t)\|_{L^2}^2+\|\nabla
\mathbf{u}_t\|_{L^2}^2 \leq C\Big(\|\nabla
\mathbf{u}\|^4_{L^2}(\|\sqrt{\rho}\mathbf{u}_t\|^2_{L^2}+\|\nabla
\mathbf{u}\|_{H^1}^2)+\|\mathbf{u}\|_{L^2}^2\Big),
\end{equation*}
which, by integrating over $(\tau,t)$ and using (\ref{0421}), leads to
\begin{equation}\label{0428}
\|\sqrt{\rho}\mathbf{u}_{t}(t)\|_{L^2}^2+\int_{\tau}^t\|\nabla
\mathbf{u}_t(s)\|_{L^2}^2\mathrm{d}s\leq
\|\sqrt{\rho}\mathbf{u}_{t}(\tau)\|_{L^2}^2 +C(T)\delta^2_0.
\end{equation}

On the other hand, multiplying (\ref{0105})$_2$ by $\mathbf{u}_t$
in $L^2$ and recalling $\mathrm{div}\,\mathbf{u}_t=0$, we find that
\begin{equation*}\label{}\begin{aligned}
\int_{\mathbb{R}^3}\rho|\mathbf{u}_t(t)|^2\mathrm{d}\mathbf{x}=
&\int_{\mathbb{R}^3}(-\varrho{g}e_3-\rho\mathbf{u}\cdot
\nabla \mathbf{u}+\mu \Delta \mathbf{u}-\nabla q)\cdot
\mathbf{u}_t\mathrm{d}\mathbf{x}\\
=&\int_{\mathbb{R}^3}(-\varrho{g}e_3-\rho\mathbf{u}\cdot \nabla
\mathbf{u}+\mu \Delta \mathbf{u})\cdot
\mathbf{u}_t\mathrm{d}\mathbf{x}, \end{aligned}
\end{equation*}
whence,
\begin{equation*}\label{0430}\begin{aligned}
\int_{\mathbb{R}^3}\rho|\mathbf{u}_t(t)|^2\mathrm{d}\mathbf{x}\leq
C\int_{\mathbb{R}^3} (\varrho^2+|\mathbf{u}|^2|\nabla
\mathbf{u}|^2+|\Delta\mathbf{u}|^2)(t)\mathrm{d}\mathbf{x}.\end{aligned}
\end{equation*}
Taking $t\to 0$ in the above inequality and using (\ref{0401}), one gets
\begin{equation*}\label{}\begin{aligned}
\limsup_{t\to 0}\int_{\mathbb{R}^3}\rho|\mathbf{u}_t(t)|^2\mathrm{d}\mathbf{x}
\leq C\delta_0^2.\end{aligned}
\end{equation*}
Therefore, letting $\tau\rightarrow 0$ in (\ref{0428}), we conclude that
\begin{equation}\label{0432}\begin{aligned}
\|\sqrt{\rho}\mathbf{u}_{t}(t)\|_{L^2}^2+\int_{0}^t\|\nabla
\mathbf{u}_t(s)\|_{L^2}^2\mathrm{d}s\leq
C(T)\delta_0^2.\end{aligned}
\end{equation}

To derive estimates of higher derivatives, we recall again that the
pair ($\mathbf{u},q$) solves the Stokes equations:
\begin{equation*}\label{} -\mu \Delta\mathbf{u} +\nabla q
=-\rho \mathbf{u}_t-\rho (\mathbf{u}\cdot \nabla \mathbf{u})-\varrho
{g}e_3,\quad \mathrm{div}\mathbf{u}=0 \mbox{ in
}\mathbb{R}^3.\end{equation*}
It follows from the classical regularity theory for Stokes equations that
\begin{equation}\label{0433}
\begin{aligned}\|\nabla^2 \mathbf{u} \|_{L^2}^2+\|\nabla q\|_{L^2}^2
\leq &C\|-\rho \mathbf{u}_t-\rho(\mathbf{u}\cdot\nabla \mathbf{u})-\varrho {g}e_3\|_{L^2}^2\\
\leq & C (\|\rho \mathbf{u}_t\|_{L^2}^2+\|\rho\mathbf{u}\cdot\nabla
\mathbf{u}\|_{L^2}^2+\|\varrho\|_{L^2}^2)\\
\leq & C(\|\rho\|_{L^\infty}
\|\sqrt{\rho}\mathbf{u}_t\|_{L^2}^2+\|\rho\|_{L^\infty}^2\|\mathbf{u}\|_{L^6}^2\|\nabla
\mathbf{u}\|_{L^3}^2+g\|\varrho\|_{L^2}^2)\\
\leq &C(\|\sqrt{\rho}\mathbf{u}_t\|_{L^2}^2+\|\nabla
\mathbf{u}\|_{L^2}^6+\|\varrho\|_{L^2}^2)+\frac{1}{2}\|\nabla
\mathbf{u}\|_{H^1}^2.
\end{aligned}\end{equation}
Putting the estimate (\ref{0421}), (\ref{0432}) and (\ref{0433}) together, we obtain
\begin{equation}\label{0436}
\|\nabla^2 \mathbf{u}(t) \|_{L^2}^2+\|\nabla q(t)\|_{L^2}^2\leq
C(T)\delta_0^2\quad\mbox{ for any }t\in (0,T].
\end{equation}

Arguing analogously to (\ref{0403}), we can deduce from (\ref{0108})$_1$ that
\begin{equation}\label{0437}
\frac{d}{dt}\|\varrho (t) \|_{L^6}\leq 6\|\bar{\rho}'\|_{L^\infty}\|\mathbf{u} (t) \|_{L^6}
\quad\mbox{ for any }t\in (0,T].
\end{equation}
By virtue of the estimates (\ref{0401}), (\ref{0421}), (\ref{0432}),
(\ref{0436}) and (\ref{0437}), and the fact that
\begin{equation*}\label{}
\begin{aligned}\|\nabla^2 \mathbf{u} \|_{L^6}^2+\|\nabla q\|_{L^6}^2\leq &C(\|\rho \mathbf{u}_t\|_{L^6}^2
+ \|\rho\mathbf{u}\cdot\nabla \mathbf{u}\|_{L^6}^2 +{g}\|\varrho
\|_{L^6}^2),
\end{aligned}\end{equation*}
we deduce that
\begin{equation}\label{0438}
\int_0^T(\|\nabla \mathbf{u}\|_{W^{1,6}}^2+\|\nabla
q\|_{L^6}^2+\|\nabla \mathbf{u}\|_{L^\infty}^2)(s)\mathrm{d}s\leq
C(T)\delta_0^2.
\end{equation}

\subsection{Estimates for $\|\varrho\|_{H^1}$}
Observing that each $\varrho$ satisfies
\begin{equation*}\label{}
(|\varrho_{x_j}|^2)_t+\mathrm{div}(|\varrho_{x_j}|^2\mathbf{u})=-2\varrho_{x_j}(
(\nabla \bar{\rho}\cdot\mathbf{u})_{x_j}+\nabla {\varrho}\cdot
\mathbf{u}_{x_j}),
\end{equation*}
we integrate over $\mathbb{R}^3$, sum over $j$ and use (\ref{0102}) to infer that
\begin{equation*}\label{}\begin{aligned}
\frac{d}{dt}\int_{\mathbb{R}^3}|\nabla
\varrho(t)|^2\mathrm{d}\mathbf{x}\leq &
C\int_{\mathbb{R}^3}\left(|\nabla \mathbf{u}||\nabla
\varrho|^2+(|\mathbf{u}|+|\nabla \mathbf{u}|\right)|\nabla
\varrho|\mathrm{d}\mathbf{x}\\
\leq &C(\|\nabla \mathbf{u}\|_{L^\infty}\|\nabla
\varrho\|_{L^2}+\|\mathbf{u}\|_{H^1})\|\nabla\varrho\|_{L^2}.\end{aligned}
\end{equation*}
Hence,
\begin{equation}\label{0442}
\|\nabla \varrho(t)\|_{L^2}\leq e^{C\int_0^t\|\nabla
\mathbf{u}(s)\|_{L^\infty}\mathrm{d}s}\left(\|\nabla
\varrho_0\|_{L^2}+C\int_0^t\|\mathbf{u}(s)\|_{H^1}\mathrm{d}s\right).
\end{equation}

From (\ref{0408}), (\ref{0438}) and (\ref{0442}) it follows that
\begin{equation*}
\| \varrho (t)\|_{H^1}\leq C(T)\delta_0\quad\mbox{ for any }t\in (0,T].
\end{equation*}
Summing up the above estimates, we arrive at the following property:
\begin{pro} \label{pro:0401}There exists a $\delta_1\in (0,1]$, such that
if $ \sqrt{\|\varrho_0\|_{H^1}^2+\|\mathbf{u}_0\|_{H^2}^2}=
\delta_0\in (0,\delta_1)$,  then any classical solution $(\varrho,
\mathbf{u}, q)$ to (\ref{0105}), emanating from the initial data
$(\varrho_0, \mathbf{u}_0)$, satisfies
\begin{equation*}\begin{aligned}\label{0444}&\sup_{0<s\leq T_\delta}(\|\varrho\|_{H^1}^2
+\|\mathbf{u}\|_{H^2}^2+\|\mathbf{u}_{t}\|_{L^2}^2+
\|\nabla q\|_{L^2}^2)(t)  \\
&+\int_0^{T_\delta}(\|\nabla \mathbf{u}\|_{H^1}^2+\|\nabla
\mathbf{u}_t\|_{L^2}^2)(s)\mm{d}s\leq C(T_\delta)\delta_0^2,
\end{aligned}\end{equation*}
where $T_\delta=\min\{T^{\max},{-{C_1}^{-1}\mathrm{ln}\delta_0}\}$
and $T^{\max}$ denotes the maximal time of existence of the
classical solution $(\varrho,\mathbf{u},q)$.
\end{pro}
\begin{rem}\label{rem:0401}
The local existence of classical solutions and global existence of
classical small solutions to the 3D nonhomogeneous incompressible Navier-Stokes
equations have been established by many authors, see \cite{PZGSS,TLYHCJ}
for example. In particular, by a slight modification, one can follow
the proof of \cite[Theorem 1.1]{TLYHCJ} and use the expanding domain technique in
 \cite{CYKHOM,CHJKH} to obtain a local existence
result of classical solutions to the Rayleigh-Taylor instability problem
(\ref{0105})--(\ref{0107}) defined on $(0,T^{\max})\times\mathbb{R}^3$,
where $T^{\max}\to\infty $ as the initial data
$\|\varrho_0\|_{H^3(\mathbb{R}^3)}+\|\mathbf{u}_0\|_{H^4(\mathbb{R}^3)}\rightarrow 0$.
The proof is standard by means of energy estimates, and hence we omit it here.
\end{rem}

\section{Proof of Theorem \ref{thm:0101}}
In this section we start to show Theorem \ref{thm:0101}. To this
end, we first construct a linear solution, and a family of nonlinear
solutions which posses some special properties. Then, by
contradiction argument, we show that there exists a nonlinearly
unstable solution satisfying the properties as stated in Theorem
\ref{thm:0101}. Suppose that $s\geq 2$, $\delta>0$, $K>0$, and $F$
satisfying (\ref{0109}), are arbitrary but given.

\subsection{Construction of a solution to the linearized problem}
 In view of Theorem \ref{thm:0203}, we can construct a
classical solution $(\varrho,\mathbf{u},q)$ to
(\ref{0106})--(\ref{0108}) satisfies the properties
(\ref{0263})--(\ref{0265}), (\ref{0268}) and (\ref{0269}). Noticing
that the solution $(\varrho,\mathbf{u})$ is independent of $s$ and
$\|{u}_3(0)\|_{H^s}>0$, we define\quad
\begin{equation}\label{0501}\begin{aligned}
(\tilde{\rho},\tilde{\mathbf{v}},\tilde{p})=
\frac{\delta({\varrho},{\mathbf{u}},q)}{\|(\varrho,\mathbf{u})(0)\|_{H^s}},
\end{aligned}\end{equation}
and find that $(\tilde{\rho},\tilde{\mathbf{v}},\tilde{p})$ is still a classical
solution to (\ref{0106})--(\ref{0108}) satisfying all the properties
of $(\varrho,\mathbf{u},q)$. Moreover,
\begin{equation*}
\|(\tilde{\rho},\tilde{\mathbf{v}})(0)\|_{H^s}=\delta.
\end{equation*}

Now, we denote
\begin{equation*}\label{}\begin{aligned}
&i_0:=i_0(s):=\frac{{\|{u}}_3(0)\|_{L^2}}{\|(\varrho,\mathbf{u})(0)\|_{H^s}}=
\frac{\|\tilde{{v}}_3(0)\|_{L^2}}{\delta}\leq {1},
\end{aligned}\end{equation*}then $i_0>0$.
Consequently, defining
$t_K=\frac{2}{\Lambda}\mathrm{ln}\frac{2K}{i_0}$ and recalling that
$\tilde{{v}}_3$ satisfies (\ref{0264}), we obtain
\begin{eqnarray}\label{0504}
\|\tilde{{v}}_3(t_K)\|_{L^2}\geq e^{t_K\Lambda/2}i_0\delta\geq 2K\delta.
\end{eqnarray}

\subsection{Construction of a solution to the corresponding nonlinear problem}
Based on the initial data $(\tilde{\rho},\tilde{\mathbf{v}})(0)$ of
the solution $(\tilde{\rho},\tilde{\mathbf{v}},\tilde{p})$ given in
(\ref{0501}), we proceed to construct a family of solutions to the
perturbed nonlinear problem. Define
\begin{equation*}(\bar{\varrho}_0^\varepsilon,\bar{\mathbf{u}}_0^\varepsilon):=
\varepsilon(\tilde{\rho},\tilde{\mathbf{v}})(0)\quad\mbox{ for }\varepsilon\in (0,1).
\end{equation*}

Noticing that
\begin{equation}\label{0506}(\bar{\varrho}_0^\varepsilon,\bar{\mathbf{u}}_0^\varepsilon)\in
(H^{\infty})^4,\mbox{ and } \|(\bar{\varrho}_0^{\varepsilon},
\bar{\mathbf{u}}_0^{\varepsilon})\|_{H^s}<
\delta\varepsilon<\delta,\end{equation} we see by Remark
\ref{rem:0401} and Proposition \ref{pro:0401} that there exists a
constant $\varepsilon_1$, such that for any $\varepsilon\in
(0,\varepsilon_1)$,
 there exists a classical solution
 $(\bar{\varrho}^\varepsilon,\bar{\mathbf{u}}^\varepsilon,\bar{q}^\varepsilon)$
to the nonlinear RT problem (\ref{0105})--(\ref{0107}) on
$(0,T_\varepsilon^{\max})\times \mathbb{R}^3$ with
$T_\varepsilon^{\max}>t_K$, satisfying
\begin{equation}\begin{aligned}\label{0507}&\sup_{0<t\leq T_K}
(\|\bar{\varrho}^\varepsilon(t)\|_{H^1}^2
+\|\bar{\mathbf{u}}^\varepsilon(t)\|_{H^2}^2+\|\bar{\mathbf{u}}_{t}^\varepsilon(t)\|_{L^2}^2
+\|\nabla \bar{q}^\varepsilon(t)\|_{L^2}^2)\\
&+\int_0^{T_K}(\|\nabla\bar{\mathbf{u}}^\varepsilon(s)\|_{H^1}^2+\|\nabla
\bar{\mathbf{u}}_t^\varepsilon(s)\|_{L^2}^2)\mm{d}s\leq C(T_K)\delta^2\varepsilon^2,
\end{aligned}\end{equation}
where $C(T_K)$ is independent of $\varepsilon$. In addition,
\begin{equation}\begin{aligned}\label{0508}
&\sup_{0<t\leq T_K}\|(\bar{\varrho}^\varepsilon+\bar{\rho})(t)\|_{L^\infty}\leq
C(\delta)\quad\mbox{ for some constatnt } C(\delta) \mbox{ depneding }\delta,\\
&\sup_{0<t\leq T_K} \|\bar{q}^\varepsilon\|_{H^1(\Omega')}^2\leq
C(T_K,\Omega')\delta^2\varepsilon^2\quad\mbox{ for any }\Omega'\subset\subset\mathbb{R}^3.
\end{aligned}\end{equation}

Obviously, to complete the proof of Theorem \ref{thm:0101}, it suffices to
show the following lemma.
\begin{lem} \label{lem5.1}
There exists an $\varepsilon_0\in(0,\varepsilon_1)$, such that
the classical solution $(\bar{\varrho}^{\varepsilon_0},\bar{\mathbf{u}}^{\varepsilon_0})$
of (\ref{0105})--(\ref{0107}), emanating from the initial data
$(\bar{\varrho}_0^{\varepsilon_0},\bar{\mathbf{u}}_0^{\varepsilon_0})$, satisfies
\begin{equation*}
\|\bar{{u}}_3^{\varepsilon_0}(t_K)\|_{L^2}>
F(\|(\bar{\varrho}_0^{\varepsilon_0},\bar{\mathbf{u}}_0^{\varepsilon_0})\|_{H^s})
\quad\mbox{for some
}t_K\in\left(0,\frac{2}{\Lambda}{\mathrm{ln}\frac{2K}{i_0}}\right]\subset
(0,T^{\max}_{\varepsilon_0}),\end{equation*} where
$T^{\max}_{\varepsilon_0}$ denotes the maximum time of existence to
the solution
$(\bar{\varrho}^{\varepsilon_0},\bar{\mathbf{u}}^{\varepsilon_0})$.
\end{lem}
\begin{pf} We shall show the lemma by contradiction. Suppose that for any
$\varepsilon\in (0,\varepsilon_1)$, the classical solution
$(\bar{\varrho}^{\varepsilon},\bar{\mathbf{u}}^{\varepsilon},\bar{q}^{\varepsilon})$,
emanating from the initial data $(\bar{\varrho}_0^{\varepsilon},\bar{\mathbf{u}}_0^{\varepsilon})$,
satisfies\quad
\begin{equation*}   \|\bar{{u}}_3^\varepsilon(t)\|_{L^2}\leq
F(\|(\bar{\varrho}_0^\varepsilon,\bar{\mathbf{u}}_0^\varepsilon)\|_{H^s})
\quad\mbox{for any }t\in \left(0,t_K\right]\subset (0,T^{\max}_\varepsilon),
\end{equation*}
which, together with  (\ref{0506}) and (\ref{0109}), yields
\begin{equation}\label{n0506}
\|\bar{{u}}_3^\varepsilon(t)\|_{L^2}\leq
K\|(\bar{\varrho}_0^\varepsilon,\bar{\mathbf{v}}_0^\varepsilon)\|_{H^s}\leq
K\delta\varepsilon ,\quad\;\;\forall\, t\in
(0,T_\varepsilon^{\max}).\end{equation}

 We denote
 $(\tilde{{\varrho}}^\varepsilon,\tilde{{\mathbf{u}}}^\varepsilon,\tilde{{q}}^\varepsilon)
 :=(\bar{\varrho}^\varepsilon,\bar{\mathbf{u}}^\varepsilon,\bar{q}^\varepsilon)/\varepsilon$,
 then they satisfy
\begin{equation}\label{n0508}\left\{\begin{array}{l}
\tilde{\varrho}^\varepsilon_t+\tilde{\bf
u}^\varepsilon\cdot\nabla (\varepsilon\tilde{\varrho}^\varepsilon+\bar{\rho})=0,\\
(\varepsilon\tilde{\varrho}^\varepsilon+\bar{\rho})\tilde{\bf u}^\varepsilon_t
+\varepsilon(\varepsilon\tilde{\varrho}^\varepsilon+\bar{\rho})\tilde{\bf
u}^\varepsilon\cdot\nabla \tilde{\bf u}^\varepsilon+\nabla
\tilde{q}^\varepsilon+g \tilde{\varrho}^\varepsilon
e_3=\mu\Delta\tilde{\bf u}^\varepsilon,\\
\mathrm{div}\tilde{\mathbf{u}}^\varepsilon=0.\end{array}\right.\end{equation}
with initial data
\begin{equation}(\tilde{{\varrho}}^\varepsilon,\tilde{{\mathbf{u}}}^\varepsilon)(0)
=({\tilde{\rho}},{\tilde{\mathbf{v}}})(0).\end{equation}
Moreover, the following estimates hold because of (\ref{0507})--(\ref{n0506}).
\begin{eqnarray}\label{0509}&&\sup_{0<t\leq
T_K}\|(\varepsilon\tilde{\varrho}^\varepsilon+\bar{\rho})(t)\|_{L^\infty}\leq C(\delta),\quad
  \sup_{0<t\leq T_K} \|\tilde{{u}}_3^\varepsilon(t)\|_{L^2}\leq K\delta,\\
&&\sup_{0<t\leq T_K} \|\tilde{q}^\varepsilon\|_{H^1(\Omega')}^2\leq
C(T_K,\Omega')\delta^2\quad\mbox{ for any }\Omega'\subset\subset \mathbb{R}^3, \\[2mm]
&&\sup_{0<t\leq T_K}(\|\tilde{\varrho}^\varepsilon\|_{H^1}^2
+\|\tilde{\mathbf{u}}^\varepsilon\|_{H^2}^2+\|\tilde{\mathbf{u}}_{t}^\varepsilon
\|_{L^2}^2+ \|\nabla \tilde{q}^\varepsilon\|_{L^2}^2)(t)\nonumber \\
&& \quad +\int_0^{T_K}(\|\nabla\tilde{\mathbf{u}}^\varepsilon\|_{H^1}^2+\|\nabla
\tilde{\mathbf{u}}_t^\varepsilon\|_{L^2}^2)(s)\mm{d}s\leq
C(T_K)\delta^2. \label{0511}  \end{eqnarray}
The continuity equation (\ref{0108})$_1$ combined with (\ref{0509})--(\ref{0511}) immediately implies
\begin{equation}\begin{aligned}\label{0512}&\sup_{0<t\leq T_K}
\|\tilde{\varrho}^\varepsilon_t\|_{L^2}^2\leq
C(T_K)\delta^2.
\end{aligned}\end{equation}

Thus, from (\ref{0509})--(\ref{0512}) we immediately infer that there exists a subsequence
(not relabeled) of $\{(\bar{\varrho}^{\varepsilon},\bar{\mathbf{u}}^{\varepsilon},\bar{q}^{\varepsilon})\}$,
such that
\begin{equation*}\begin{aligned} &
(\tilde{\varrho}^\varepsilon_t,\tilde{\mathbf{u}}_{t}^\varepsilon,
\tilde{q}^\varepsilon)\rightarrow (\tilde{\varrho}_t,\tilde{\mathbf{u}}_{t}, \tilde{q})
\mbox{ weakly-star in }L^\infty (0,T_K;(L^2)^4\times H^1_{\mathrm{loc}}),\\
&(\tilde{\varrho}^\varepsilon,\tilde{\mathbf{u}}^\varepsilon)\rightarrow
(\tilde{\varrho},\tilde{\mathbf{u}})\mbox{ weakly-star in }L^\infty (0,T_K; H^1)\times (H^2)^3),\\
&(\tilde{\varrho}^\varepsilon,\tilde{\mathbf{u}}^\varepsilon)\rightarrow
(\tilde{\varrho},\tilde{\mathbf{u}})\mbox{ strongly in }C^0([0,T_K],(L^2_{\mathrm{loc}})^4),
\end{aligned}\end{equation*}
and
\begin{equation}\begin{aligned}\label{0519}\sup_{0<t\leq T_K}\|\tilde{{{u}}}_3(t)\|_{L^2}
\leq K\delta,\qquad (\tilde{\varrho},\tilde{\mathbf{u}})\in C^0([0,T_K],(L^2)^4).
\end{aligned}\end{equation}

If one takes to the limit as $\varepsilon\to 0$ in the equations (\ref{n0508}), one gets
\begin{equation*}\label{}\left\{\begin{array}{l}
\partial_t\tilde{\varrho}+\tilde{\bf
u}\cdot\nabla \bar{\rho}=0,\\
\bar{\rho}\partial_t \tilde{\bf u}+\nabla \tilde{q}+g
\tilde{\varrho} e_3=\mu\Delta\tilde{\bf u},\\
\mathrm{div}\tilde{\mathbf{u}}=0.\end{array}\right.\end{equation*}
Therefore, we see that $(\tilde{\varrho},\tilde{\mathbf{u}})$ is just a strong
solution of the linearized problem (\ref{0106})--(\ref{0108}). Of
course, $(\tilde{\rho},\tilde{\mathbf{v}})$ is also a strong
solution of (\ref{0106})--(\ref{0108}). Moreover, $(\tilde{\varrho},
\tilde{\mathbf{u}})(0)=(\tilde{\rho},\tilde{\mathbf{v}})(0)$. Hence, according to Theorem \ref{thm:0301},
\begin{equation*}\label{0566}
\tilde{{\mathbf{v}}}=\tilde{\mathbf{u}}\mbox{ on }[0,T_K]\times
\mathbb{R}^3.
\end{equation*}
Thus, we may chain together (\ref{0504}) and the inequality (\ref{0519}) to get
\begin{equation*}\label{0567}
2K\delta\leq \| \tilde{{v}}_3(t_K)\|_{L^2(\mathbb{R}^3)}\leq \|
\tilde{{u}}_3(t_K)\|_{L^2(\mathbb{R}^3)}\leq K\delta,
\end{equation*}
which is a contraction. This completes the proof of Lemma \ref{lem5.1}, and hence the proof
of Theorem \ref{thm:0101}.
\hfill $\Box$
\end{pf}

\renewcommand\refname{References}
\renewenvironment{thebibliography}[1]{%
\section*{\refname}
\list{{\arabic{enumi}}}{\def\makelabel##1{\hss{##1}}\topsep=0mm
\parsep=0mm
\partopsep=0mm\itemsep=0mm
\labelsep=1ex\itemindent=0mm
\settowidth\labelwidth{\small[#1]}%
\leftmargin\labelwidth \advance\leftmargin\labelsep
\advance\leftmargin -\itemindent
\usecounter{enumi}}\small
\def\newblock{\ }
\sloppy\clubpenalty4000\widowpenalty4000
\sfcode`\.=1000\relax}{\endlist}
\bibliographystyle{model1-num-names}

\begin{thebibliography}{21}
\expandafter\ifx\csname
natexlab\endcsname\relax\def\natexlab#1{#1}\fi
\providecommand{\bibinfo}[2]{#2} \ifx\xfnm\relax
\def\xfnm[#1]{\unskip,\space#1}\fi
\bibitem[{Chandrasekhar(1961)}]{CSHHS}
\bibinfo{author}{S.~Chandrasekhar}, \bibinfo{title}{{Hydrodynamic and
  Hydromagnetic Stability, The International Series of Monographs on Physics}},
  \bibinfo{publisher}{Oxford, Clarendon Press}, \bibinfo{year}{1961}.
\bibitem[{Choe and Kim(2003)}]{CHJKH}
\bibinfo{author}{H.J. Choe}, \bibinfo{author}{H.~Kim},
\newblock \bibinfo{title}{Strong solutions of the Navier-Stokes equations for
  nonhomogeneous incompressible fluids},
\newblock \bibinfo{journal}{Comm. PDE}
  \bibinfo{volume}{28} (\bibinfo{year}{2003}) \bibinfo{pages}{1183--1201}.
\bibitem[{Cho and Kim(2006)}]{CYKHOM}
\bibinfo{author}{Y.~Cho}, \bibinfo{author}{H.~Kim},
\newblock \bibinfo{title}{On classical solutions of the compressible
  Navier-Stokes equations with nonnegative initial densities},
\newblock \bibinfo{journal}{Manuscripta Math.} \bibinfo{volume}{120}
  (\bibinfo{year}{2006}) \bibinfo{pages}{91--129}.
\bibitem[{Dautray and Lions(1985)}]{RDJLLA}
\bibinfo{author}{R.~Dautray}, \bibinfo{author}{J.L. Lions},
  \bibinfo{title}{Analyse math\'ematique et calcul mum\'erique pour les sciences et les techniques},
  \bibinfo{address}{Masson, Paris}, \bibinfo{year}{1985}.

\bibitem[{Duan et~al.(2011)Duan, Jiang, and Jiang}]{RDFJSJO}
\bibinfo{author}{R.~Duan}, \bibinfo{author}{F.~Jiang},\bibinfo{author}{~S.~Jiang},
\newblock \bibinfo{title}{On the Rayleigh-Taylor instability for incompressible
magnetohydrodynamic flows},
\newblock \bibinfo{journal}{SIAM J. Appl. Math.} \bibinfo{volume}{71}
  (\bibinfo{year}{2011}) \bibinfo{pages}{1990-2013}.

\bibitem[{Duan et~al.(2011)Duan, Jiang, and Jiang}]{RDFJSJS}
\bibinfo{author}{R.~Duan}, \bibinfo{author}{F.~Jiang},\bibinfo{author}{~S.~Jiang},
\newblock \bibinfo{title}{Rayleigh-Taylor instability for compressible rotating flows},
\newblock \bibinfo{journal}{submitted}
  \bibinfo{year}{2011}.

\bibitem[{Erban(1987)}]{EDGTC1111}
\bibinfo{author}{D.~Erban},
\newblock \bibinfo{title}{The equations of motion of a perfect fluid with free
  boundary are not well posed},
\newblock \bibinfo{journal}{Comm. PDE}
  \bibinfo{volume}{12} (\bibinfo{year}{1987}) \bibinfo{pages}{1175--1201}.
\bibitem[{Galdi(1994)}]{GGPA94}
\bibinfo{author}{G.~Galdi}, \bibinfo{title}{{An introduction to the
  mathematical theory of the Navier-Stokes equations. Linearized Steady
  Problems. Springer Tracts in Natural Philosophy 38. Vol. 1.}},
  \bibinfo{publisher}{Springer-Verlag}, \bibinfo{address}{New York},
  \bibinfo{year}{1994}.
\bibitem[{Guo and Tice(2011{\natexlab{a}})}]{GYTI1}
\bibinfo{author}{Y.~Guo}, \bibinfo{author}{I.~Tice},
\newblock \bibinfo{title}{Compressible, inviscid Rayleigh-Taylor instability},
\newblock \bibinfo{journal}{arXiv:0911.4098v2 [math.AP] 23 Feb 2011, to appear
  in Indiana Univ. Math. J.}  (\bibinfo{year}{2011}{\natexlab{a}}).
\bibitem[{Guo and Tice(2011{\natexlab{b}})}]{GYTI2}
\bibinfo{author}{Y.~Guo}, \bibinfo{author}{I.~Tice},
\newblock \bibinfo{title}{Linear Rayleigh-Taylor instability for viscous,
  compressible fluids},
\newblock \bibinfo{journal}{SIAM J. Math. Anal.} \bibinfo{volume}{42}
  (\bibinfo{year}{2011}{\natexlab{b}}) \bibinfo{pages}{1688--1720}.

\bibitem[{Hwang(2008)}]{HHVQ}
\bibinfo{author}{H.J. Hwang},
\newblock \bibinfo{title}{{Variational approach to nonlinear gravity-driven
  instability in a MHD setting}},
\newblock \bibinfo{journal}{Quart. Appl. Math.} \bibinfo{volume}{66}
  (\bibinfo{year}{2008}) \bibinfo{pages}{303--324}.
\bibitem[{Hwang and Guo(2003)}]{HHJGY}
\bibinfo{author}{H.J. Hwang}, \bibinfo{author}{Y.~Guo},
\newblock \bibinfo{title}{{On the dynamical Rayleigh-Taylor instability}},
\newblock \bibinfo{journal}{Arch. Rational Mech. Anal.} \bibinfo{volume}{167}
  (\bibinfo{year}{2003}) \bibinfo{pages}{235--253}.
\bibitem[{Jang and Tice(2011)}]{JJHTI}
\bibinfo{author}{J.~Jang}, \bibinfo{author}{I.~Tice},
\newblock \bibinfo{title}{Instability theory of the Navier-Stokes-Poisson equations},
\newblock \bibinfo{journal}{arXiv:1105.5128v2 [math.AP] 13 Jun 2011}
  (\bibinfo{year}{2011}).
\bibitem[{Jiang et~al.(2011)Jiang, Jiang, and Wang}]{JFJSWWWO}
\bibinfo{author}{F.~Jiang}, \bibinfo{author}{S.~Jiang},
  \bibinfo{author}{W.~Wang},
\newblock \bibinfo{title}{{On the Rayleigh-Taylor instability for two uniform
  viscous incompressible flows}},
\newblock \bibinfo{journal}{submitted} (\bibinfo{year}{2011}).
\bibitem[{Jiang et~al.(2011)Jiang, Jiang, and Wang}]{JFJSWYJ}
\bibinfo{author}{F.~Jiang}, \bibinfo{author}{S.~Jiang},
  \bibinfo{author}{Y.~Wang},
\newblock \bibinfo{title}{{On the Rayleigh-Taylor instability for
  incompressible, inviscid magnetohydrodynamic flows}},
\newblock \bibinfo{journal}{submitted}  (\bibinfo{year}{2011}).
\bibitem[{Kruskal and Schwarzschild(1954)}]{KMSMSP}
\bibinfo{author}{M.~Kruskal}, \bibinfo{author}{M.~Schwarzschild},
\newblock \bibinfo{title}{Some instabilities of a completely ionized plasma},
\newblock \bibinfo{journal}{Proc. Roy. Soc. (London) A} \bibinfo{volume}{233}
  (\bibinfo{year}{1954}) \bibinfo{pages}{348--360}.
\bibitem[{Lions(1996)}]{LPLMTFM96}
\bibinfo{author}{P.L.~Lions}, \bibinfo{title}{{Mathematical Topics in Fluid
  Mechanics: Incompressible models}}, \bibinfo{publisher}{Oxford University
  Press}, \bibinfo{address}{USA}, \bibinfo{year}{1996}.
\bibitem[{Pr$\mathrm{\ddot{u}}$ess and Simonett(2009)}]{JPGS09}
\bibinfo{author}{J.~Pr\"{u}ess}, \bibinfo{author}{G.~Simonett},
\newblock \bibinfo{title}{{On the Rayleigh-Taylor instability for the
  two-phase Navier-Stokes equations}},
\newblock \bibinfo{journal}{arXiv: 0908.3334v1 [math.AP] 23 Aug 2009}
  (\bibinfo{year}{2009}).
\bibitem[{Rayleigh(1883)}]{RLAP}
\bibinfo{author}{L.~Rayleigh},
\newblock \bibinfo{title}{{Analytic solutions of the Rayleigh equations for
  linear density profiles}},
\newblock \bibinfo{journal}{Proc. London. Math. Soc.} \bibinfo{volume}{14}
  (\bibinfo{year}{1883}) \bibinfo{pages}{170--177}.
\bibitem[{Rayleigh(1990)}]{RLIS}
\bibinfo{author}{L.~Rayleigh},
\newblock \bibinfo{title}{{Investigation of the character of the equilibrium of
  an in compressible heavy fluid of variable density}},
\newblock \bibinfo{journal}{Scientific Paper, II}  (\bibinfo{year}{1990})
  \bibinfo{pages}{200--207}.
\bibitem[{Taylor(1950)}]{TGTP}
\bibinfo{author}{G.~I. Taylor},
\newblock \bibinfo{title}{{The stability of liquid surface when accelerated in
  a direction perpendicular to their planes}},
\newblock \bibinfo{journal}{Proc. Roy Soc. A} \bibinfo{volume}{201}
  (\bibinfo{year}{1950}) \bibinfo{pages}{192--196}.
\bibitem[{Tong and Yuan(2010)}]{TLYHCJ}
\bibinfo{author}{L.~Tong}, \bibinfo{author}{H.~Yuan},
\newblock \bibinfo{title}{{Classical solutions to Navier-Stokes equations for
  nonhomogeneous incompressible fluids with non-negative densities}},
\newblock \bibinfo{journal}{J. Math. Anal. Appl.} \bibinfo{volume}{362}
  (\bibinfo{year}{2010}) \bibinfo{pages}{476--504}.
\bibitem[{Wang(1994)}]{WJH}
\bibinfo{author}{J.~Wang}, \bibinfo{title}{Two-Dimensional Nonsteady Flows and
  Shock Waves (in Chinese)}, \bibinfo{publisher}{Science Press},
  \bibinfo{address}{Beijing, China}, \bibinfo{year}{1994}.
\bibitem[{Wang(2010)}]{WYC}
\bibinfo{author}{Y.~Wang},
\newblock \bibinfo{title}{{Critical magnetic number in the MHD Rayleigh-Taylor
  instability}},
\newblock \bibinfo{journal}{arXiv:1009.5422v1 [math.AP] 28 Sep 2010}
  (\bibinfo{year}{2010}).


\bibitem[{Zhang(2008)}]{PZGSS}
\bibinfo{author}{P.~Zhang},
\newblock \bibinfo{title}{{Global smooth solutions to the 2D nonhomogeneous
  Navier-Stokes equations}},
\newblock \bibinfo{journal}{Inter. Math. Research Notices}
  \bibinfo{volume}{Vol. 2008, Article ID rnn098, 26 pages.
  doi:10.1093/imrn/rnn098} (\bibinfo{year}{2008}).

\end{thebibliography}

\end{document}